\documentclass[11pt]{article}

\usepackage{graphicx}
\usepackage{amssymb}
\usepackage{amsmath}
\usepackage{amsthm}
\usepackage{latexsym}
\usepackage[english]{babel}
\usepackage{tikz}
\usetikzlibrary{matrix,arrows,decorations.pathmorphing}
\usepackage{tikz-cd}
\usepackage{appendix}

\newtheorem{theorem}{Theorem}[section]
\newtheorem{definition}{Definition}[section]
\newtheorem{lemma}{Lemma}[section]
\newtheorem{proposition}{Proposition}[section]
\newtheorem{remark}{Remark}[section]
\newtheorem{example}{Example}[section]
\newtheorem{corollary}{Corollary}[section]

\def \C {\mathbb{C}}
\def \R {\mathbb{R}}

\def \N {\mathbb{N}}
\def \Z {\mathbb{Z}}

\def \D {\mathcal{D}}
\def \Q {\mathcal{Q}}

\def\XXint#1#2#3{{\setbox0=\hbox{$#1{#2#3}{\int}$ }
\vcenter{\hbox{$#2#3$ }}\kern-.6\wd0}}

\begin{document}

\title{Conformal Fractional Dirac Operator and Spinorial $Q$-curvature}

\author{Ali Maalaoui$^{(1)(2)}$}
\addtocounter{footnote}{1}
\footnotetext{Department of Mathematics, Clark University, 950 Main Street, Worcester, MA 01610, USA. E-mail address: 
{\tt{amaalaoui@clarku.edu}}}
\addtocounter{footnote}{1}
\footnotetext{Department of Mathematics, MIT, 77 Massachusetts Avenue
Cambridge, MA 02139-4307. E-mail address: 
{\tt{maala650@mit.edu}}}

\date{}
\maketitle

\vspace{5mm}

{\noindent\bf Abstract} {In this paper we introduce the conformal fractional Dirac operator and its associated fractional spinorial Yamabe problem. We also present a Caffarelli-Silvestre type extension for this fractional operator, allowing us to express it as a Dirichlet-to-Neumann type operator. As a consequence, we exhibit energy inequalities associated to this operator along with a weighted type Sobolev inequality for spinors. In the second part of the paper, we focus on the critical operator (which can be local or non-local depending on the evenness of the dimension). We introduce a $Q$-curvature operator, acting on spinors generalizing the classical notion of the scalar $Q$-curvature.}

\vspace{5mm}

\noindent
{\small Keywords: Fractional Operator, Dirac Operator, Q-curvature}

\vspace{4mm}

\noindent
{\small 2010 MSC. Primary: 53C27, 30F45.  Secondary: 58J32, 53C18.}

\vspace{4mm}


\section{Introduction}

In the recent years there have been a great interest in fractional operators in conformal and CR geometries and their related scalar invariants. For instance, we refer to the work on the conformal fractional Laplacian and the fractional $Q$-curvature in \cite{Gon,CC,CG, MQ} and the conformal fractional sub-Laplacian in \cite{RGM}. The primary interest of this paper is to define and study the conformal fractional Dirac operator on spin manifolds. Indeed, in \cite{GMP2, F, F2,FKS} there have been studies and constructions of conformally invariant powers of the Dirac operator in a way similar to the classical GJMS operators (\cite{GJMS}). But unlike the scalar case, the fractional setting in spinorial geometry is still not well understood. We recall that the fractional GJMS operators consist of a family of conformally covariant pseudo-differential operators defined on the boundary of a Poicare-Einstein manifold via scattering theory as detailed in \cite{GZ} and can be formulated as a generalized Dirichlet-to-Neumann operator associated to weighted GJMS operators of suitable order defined on the interior of a conformal compactification of the Poincar\'{e}-Einstein structure. Our first objective is to define the conformal fractional Dirac operator $\D_{h}^{2\lambda}$ using a scattering operator adapted to spinors \cite{GMP,GMP2} and then formulate a Caffarelli-Silvestre type extension allowing us to see this spinorial fractional operator as a Dirichlet-to-Neumann type operator for a second order weighted operator in the interior. The study of the conformal fractional Dirac operator stems from the similarities that exist between the spinorial Yamabe problem, consisting of optimizing the first eigenvalue of the Dirac operator in a conformal class, see \cite{Amm0,Amm1,AGHM,AH} and the classical Yamabe problem which also can be seen as an optimization problem for the first eigenvalue of the conformal Laplacian in a given conformal class. It is then natural to attempt an extension of the problem to the fractional setting as it was the case for the fractional Yamabe problem (see \cite{MQ}). But, in order to be able to make effective computations with these fractional operators, the Caffarelli-Silvestre type extension, see \cite{CS}, becomes fundamental. The second question that we address deals with some energy estimates relating the energy of $\D_{h}^{2\lambda}$ and a weighted energy of the extension, focusing mainly on the case $\lambda \in (0,\frac{1}{2})$. In particular, we prove a weighted Sobolev type inequality exhibiting a surprising link between these operators and the fractional $Q$-curvature.
In the second part of this paper we focus on the critical case corresponding to $\lambda=\frac{n}{2}$, where we have different behaviour depending on the evenness of the dimension $n$. Indeed, we define a version of spinorial $Q$-curvature (in this case it is a spinor valued operator) and an associated scalar invariant. In fact, this construction is similar to the one done for differential forms on even dimensional manifolds \cite{BG,BG2, AG}. Indeed, the authors prove the existence of a sequence of conformally invariant operators on forms and their associated $Q$-curvature. A similar construction was also established for the Yang-Mills functional on $6$-dimensional Poincar\'{e}-Einstein manifolds \cite{GPS}. \\

The structure of this paper is as follows:\\
Section 2 of this paper serves as the main motivation and background that lead to the investigation of the conformal fractional Dirac operator. Indeed, in Sub-Section 2.1, we recall some results related to the classical Spinorial Yambe problem  such as the Aubin-type inequality and the solvability of the problem in the case of a non-locally conformally flat manifold of dimension greater than six. We also discuss the existence of singular solutions in the case of the punctured sphere. In Sub-Section 2.2, we survey the construction of the fractional Laplacian via the scattering operator on Poicare-Einstein manifolds, along with the GJMS operators that appear as the residue of the scattering operator at even integers of the parameter. This review is meant to guide the intuition of the reader when dealing with the spinorial case. In Sub-Section 2.3 we provide some background about spinorial geometry and the scattering operator associated to the Dirac operator on a Poincar\'{e}-Einstein manifold. In Sub-Section 2.4, we define the Conformal fractional Dirac operator $\D_{h}^{2\lambda}$ and we formulate the fractional spinorial Yamabe problem in Theorem \ref{thmyam}. We also provide two examples where we have an explicit formula for $\D_{h}^{2\lambda}$. Namely, we provide the expression of the operator on the standard sphere and the flat Euclidean space. Moreover, we provide the standard "bubble" for the fractional Spinorial Yamabe problem in these two cases.\\

In Section 3, we focus on the Caffarelli-Silvestre type extension, which allows us to formulate the conformal fractional Dirac operator as a Dirichlet-to-Neumann type operator for a boundary value problem. In fact, the main extension result is stated in Theorem \ref{thmext1}, for a general defining function. Moreover, for the case $\lambda \in (0,\frac{1}{2})$ we show in Theorem \ref{thmext2}, that using a $\lambda$-admissible defining function, the Dirichlet-to-Neumann type operator involves the fractional $Q$-curvature, $Q_{h}^{\lambda}$, of the conformal boundary at infinity. One can also formulate the extension problem using a first order operator involving the weighted Dirac operator.\\

In Section 4, we present several energy identities relating the conformal fractional Dirac operator and the extension to the compactified manifold $\overline{X}^{n+1}$. In Proposition \ref{enbas}, we state the basic energy identities that follow from the scattering definition of the operator in question. Then, in Theorem \ref{en1}, we state the optimized energy identities corresponding to the extension problem introduced in Section 3. In Sub-Section 4.1, we focus in Proposition \ref{en2}, on the energy identities that are tied to the use of the $\lambda$-admissible adapted function and the corresponding fractional $Q$-curvature. Moreover, we provide an improved identity for the case $\lambda \in (0,\frac{1}{4})$, connecting the square of the conformal fractional Dirac operator to the $Q$-curvature $Q_{h}^{2\lambda}$, in Theorem \ref{en3}. In Sub-Section 4.2, we prove a weighted Sobolev type inequality for spinors, stated in Theorem \ref{sob} and characterize the equality case.\\

In Section 5, we focus on the critical case corresponding to $\lambda=\frac{n}{2}$. There, we define a spinorial $Q$-curvature type operator and show its variation under conformal change when restricted to the kernel of the critical operator $\D_{h}^{n}$, as stated in Theorem \ref{thmQ}. Among the results that we prove is the existence of a scalar invariant related to spinors in the kernel. Moreover, in Theorem \ref{thmQ2}, we present a construction of a Spinor on the Poincar\'{e}-Einstein manifold with a specific asymptotic behavior similar to the construction of the special defining function of Graham-Fefferman \cite{FG1}. \\

{\noindent\bf Acknowledgment}

\noindent
The author is supported by the AMS-Simons Research Enhancement Grant for PUI faculty under the project "Conformally Invariant
Non-Local Equations on Spin Manifolds". He also wants to express his gratitude to the department of Mathematics at MIT for the warm hospitality during the work on this manuscript.
\section{Motivation and Preliminaries}
\subsection{The Spinorial Yamabe Problem}
Let $(M,h,\Sigma_{h}M)$ be a compact spin manifold. There is a natural differential operator acting on smooth sections of $\Sigma_{h}M$, namely, the Dirac operator $D_{h}$. We recall that locally, one can express $D_{h}$ as follows:
$$D_{h}\psi=\sum_{k=1}^{n}e_{k}\cdot \nabla_{e_{k}}\psi,$$
where $(e_{1},\cdots,e_{n})$ is an $h$-orthonormal frame. The Dirac operatoris a strong tool that is classically used in understanding the geomery and topology of a manifold.

The Dirac operator has also many applications in conformal geomery, mainly because of its conformal invariance property. That is, if $\overline{h}=e^{2u}h$ then there exists an isometric identification between $\Sigma_{h}$ and $\Sigma_{\overline{h}}$ so that
\begin{equation}\label{iden}
D_{\overline{h}}\overline{\psi}=\overline{e^{\frac{n+1}{2}u}D_{h}e^{-\frac{n-1}{2}u}\psi},
\end{equation}
In particular, $\ker D_{h}$ is a conformal invariant. In what follows, we will disregard the isometric identification and its notation.
In \cite{Amm0,AH,AGHM}, the authors formulated a version of the Yamabe problem for spinor. This formulation can be seen as an optimization problem of the first eigenvalue of $D_{h}$ under a constant volume constraint. Indeed, if $\lambda_{1}(D_{h})$ is the first non-zero eigenvalue of  $D_{h}$, then 
$$\lambda^{+}(M,[h]):=\inf_{\hat{h}\in [h]}|\lambda_{1}(D_{\hat{h}})|Vol(\hat{h})^{\frac{1}{n}}.$$
This conformal invariant is called the  B\"{a}r-Hijazi-Lott invariant and it has the following properties:
\begin{itemize}
\item[i) \cite{Amm0}] If  $Y(M,[h])$ is the Yamabe constant of $(M,[h])$, then $$(\lambda^{+}(M,[h]))^{2}\geq \frac{n}{4(n-1)}Y(M,[h]).$$
\item[ii) \cite{Amm0,AGHM}] $\lambda^{+}(M,[h])\leq \lambda^{+}(S^{n},[g_{0}])$ and if the inequality is strict, then there exists a generalized metric $\overline{h}\in [h]$ with singular set $S\subset  M$ and a spinor $\phi$ such that $|\phi|=1$ on $M\setminus S$ such that $D_{\overline{h}}\phi=\lambda^{+}\phi$.
\item[iii) \cite{SX}] If $n\geq 6$ and $(M,h)$ is not locally conformally flat then $\lambda^{+}(M,[h])< \lambda^{+}(S^{n},[g_{0}])$.
\item[iv) \cite{MSX}] There exists a one parameter family of metrics $(g_{T})_{T>T_{0}}$ on $S^{1}\times S^{n-1}$ such that for $T\to +\infty$
$$\lambda^{+}(S^{1}\times S^{n-1},[g_{T}])\to \lambda^{+}(S^{n},[g_{0}])$$
\end{itemize}
These properties show in particular the deep similarities between the Spinorial Yamabe problem and the classical Yamabe problem. As we know, the Yamabe problem has higher order and fractional order extensions. Our objective is to state these extensions in the spinorial setting. For that purpose, we will review the construction of the fractional conformal Laplacian and then move to the construction of the fractional Dirac operator.\\

\subsection{The scattering operator and fractional Laplacian}
We recall that Poincar\'{e}-Einstein manifold is a triple $(X^{n+1},M^{n},g_{+})$  such that $(X^{n+1},g_{+})$ is an Einstein manifold with $Ric_{g_{+}}=-ng_{+}$ and $X^{n+1}$ is diffeomeorphic to the interior of a compact manifold $\overline{X}^{n+1}$ such that $\partial \overline{X}^{n+1}=M$. Moreover, we assume the existence of a a smooth positive defining function $x:\overline{X}^{n+1}\to \R_{+}$ such that $M=x^{-1}(0)$ and $\overline{g}:=x^{2}g_{+}$ extends to a $C^{n,\alpha}$ metric on $\overline{X}^{n+1}$ with $|dx|_{\overline{g}}=1$ on $M$. If $w \in C^{\infty}(\overline{X}^{n+1})$ then $e^{w}x$ is also another defining function. Therefore, at the boundary $M$, only the conformla class $[\overline{g}|_{TM}]$ is well defined. We say that $(M,[\overline{g}|_{TM}])$ is the conformal infinity of the Poincar\'{e}-Einstein manifold $(X^{n+1},g_{+})$. For each element $h\in [\overline{g}|_{TM}]$ corresponds a defining function $x$ on $(\overline{X}^{n+},\overline{g})$, unique in the neighborhood of $M$, called a geodesic defining function such that $g_{+}=\frac{dx^{2}+h_{x}}{x^{2}}$, where $h_{x}$ is a one parameter family of metrics on $M$ with $h_{0}=h$. The family $h_{x}$ has the following expansion near $M$:
$$\left\{\begin{array}{ll}
h_{x}=h+h_{(2)}x^{2}+\cdots +h_{(n-1)}x^{n-1}+Lx^{n} +o(x^{n}), \text{ if $n$ is odd}\\
h_{x}=h+h_{(2)}x^{2}+\cdots +h_{(n-2)}x^{n-2}+h_{(n)}x^{n}\ln(x)+Lx^{n} +o(x^{n}), \text{ if $n$ is even}.
\end{array}
\right.
$$
The $h_{(\ell)}$ for $\ell \leq n$ and even, are locally determined by the metric $h$. But the coefficient $k$ is nonlocal.

We focus now on the spectral properties the Laplace operator $-\Delta_{g_{+}}$ on $(X^{n+1},g_{+})$. Indeed, we know from \cite{GZ} that $Spec_{ess}(-\Delta_{g_{+}})\subset [\frac{n^{2}}{4},+\infty)$ and the point spectrum $Spec_{pp}(-\Delta_{g_{+}})\subset (0,\frac{n^{2}}{4})$. The study of generalized eigenfunctions of the operator $-\Delta_{g_{+}}$ is the main ingredient to the construction of conformally invariant operator on the boundary at infinity. We describe here the construction of the classical GJMS operator using the the Scattering operator as detailed in \cite{GZ}. In fact, one can see the GJMS operators as an obstruction to the regularity of generalized eigenfunctions of $-\Delta_{g_{+}}$. We will summarize the construction for the scalar case since it will be relevant for the rest of our work.
\begin{theorem}[\cite{GZ}]\label{thmscal}
Let $f\in C^{\infty}(M)$ and $\gamma \in (0,\frac{n}{2})$ then the eigenvalue problem
$$-\Delta_{g_{+}}u-(\frac{n}{2}-\gamma)(\frac{n}{2}+\gamma)u=0$$
has a unique solution of the form
$$u=\left\{\begin{array}{ll}
x^{\frac{n}{2}-\gamma}F+x^{\frac{n}{2}+\gamma}G \text{ if $\gamma \not \in \N$}\\
x^{\frac{n}{2}-\gamma}A+x^{\frac{n}{2}+\gamma}\ln(x)B \text{ if $\gamma \in \N$}
\end{array}
\right.$$
where $F, G, A, B \in C^{\infty}(\overline{X})$ and $F|_{M}=A|_{M}=f$.
\end{theorem}
Based on this theorem, the scattering operator $S(\gamma)$ is defined for $\gamma \not \in \N$ and such that $(\frac{n^{2}}{4}-\gamma^{2})\not \in Spec_{pp}(-\Delta_{g_{+}})$, by $S(\gamma):F|_{x=0}\to G|_{x=0}$.
\begin{proposition}[\cite{GZ}]
The scattering operator is finite meromorphic for $\Re(\gamma)>0$ with simple poles at $\gamma =k \in \N$ and
$$Res_{\gamma=k}S(\gamma)=c_{k}P_{k},$$
where $P_{k}$ is the GJMS operator of order $2k$. Moreover, using the notations of Theorem \ref{thmscal}, $$B|_{x=0}=c_{k}P_{k}.$$ 
\end{proposition}
We recall that for instance $P_{1}=L_{h}=-\Delta_{h}+\frac{(n-2)}{4(n-1)}R_{h}$ and $P_{2}=(-\Delta_{h})^{2}+\delta(a_{n}R_{h}h+b_{n}Ric_{h})d+\frac{n-4}{2}Q_{h}^{2}$, where $a_{n}$ and $b_{n}$ are constants depending on $n$ and $Q_{h}^{2}$ is the classical fourth order $Q$-curvature.
The fractional Laplacian $P_{h}^{\gamma}$ is then defined by
$$P_{h}^{\gamma}=c_{\gamma}S(\gamma),$$
where $c_{\gamma}=2^{2\gamma}\frac{\Gamma(\gamma)}{\Gamma(-\gamma)}$. The associated scalar invariant, is the $Q$-curvature, defined by
$$(\frac{n}{2}-\gamma)Q_{h}^{\gamma}:=P_{h}^{\gamma}(1).$$
The fractional Laplacian and the $Q$-curvature were extensively studied in the past few years for instance we refer the reader to \cite{Gon, MQ,Case1,CC,CG} and the references therein.
\begin{definition}\label{admis}
Given $\gamma \in (0,1)$, we say that $\rho$ is a $\gamma$-admissible defining function for the Poincar\'{e}-Einstein manifold $(X^{n+1},g_{+})$, if there exists $W\in C^{\infty}(M)$ such that $\rho=x\Big(1+Wx^{2\gamma}+o(x^{2\gamma})\Big)$. 
\end{definition}
This definition extends the notion of defining functions to the non-smooth category (it is smooth in the interior but not up to the boundary). Such defining functions are very important and had extensive use in \cite{Case1,MQ}. For instance, one can always construct a $\gamma$-admissible function by taking $v$ the solution of
$$\left\{\begin{array}{ll}
-\Delta_{g_{+}}v-(\frac{n}{2}-\gamma)(\frac{n}{2}+\gamma)v=0\\
x^{\gamma -\frac{n}{2}}v|_{x=0}=1.
\end{array}
\right.
$$
Then $\rho=v^{\frac{2}{n-2\gamma}}$ is a $\gamma$-admissible defining function, since $\rho=x\Big(1+\frac{Q_{h}^{\gamma}}{c_{\gamma}}x^{2\gamma}+o(x^{2\gamma})\Big)$. For this particular $\gamma$-admissible defining function, the scalar curvature of $\tilde{g}=\rho^{2}g_{+}$ has a nice explicit formula. Indeed, for a general defining function $x$, we have that
$$R_{\overline{g}}=n(n+1)\frac{(|\nabla^{\overline{g}}x|^{2}-1)}{x^{2}}-2n\frac{\Delta_{\overline{g}}x}{x}.$$
But for the metric $\tilde{g}=\rho^{2}g_{+}$ we have
\begin{equation}\label{curv}
R_{\tilde{g}}=n(1-2\lambda)\frac{(|\nabla^{\tilde{g}}\rho|^{2}-1)}{\rho^{2}}=\frac{2n(1-2\lambda)}{n+2\lambda}\frac{\Delta_{\tilde{g}}\rho}{\rho}.
\end{equation}

\subsection{Spinorial Geometry and the Scattering Operator}
Given a $n+1$-dimensional Riemannian manifold (possibly with boundary) $(X,g)$, we recall that a spin structure is a pair $(P_{Spin}(X,g),\sigma)$, where $P_{Spin}(X,g)$ is a $Spin(n+1)$-principal bundle and $\beta : P_{Spin}(X,g)\to P_{SO}(X,g)$ is a 2-fold covering map, which restricts to a non-trivial covering $\kappa: Spin(n+1)\to SO(n+1)$ on each fiber. Therefore, the quotient of each fiber by $\Z_{2}$ is isomorphic to the frame bundle of $X$ and hence, the following diagram commutes:
\begin{center}

\tikzset{every picture/.style={line width=0.75pt}} 

\begin{tikzpicture}[x=0.75pt,y=0.75pt,yscale=-1,xscale=1]

\draw    (173,59) -- (380,59) ;
\draw [shift={(382,59)}, rotate = 180] [color={rgb, 255:red, 0; green, 0; blue, 0 }  ][line width=0.75]    (10.93,-3.29) .. controls (6.95,-1.4) and (3.31,-0.3) .. (0,0) .. controls (3.31,0.3) and (6.95,1.4) .. (10.93,3.29)   ;
\draw    (150,80) -- (247.36,147.86) ;
\draw [shift={(249,149)}, rotate = 214.88] [color={rgb, 255:red, 0; green, 0; blue, 0 }  ][line width=0.75]    (10.93,-3.29) .. controls (6.95,-1.4) and (3.31,-0.3) .. (0,0) .. controls (3.31,0.3) and (6.95,1.4) .. (10.93,3.29)   ;
\draw    (415,79) -- (314.63,149.85) ;
\draw [shift={(313,151)}, rotate = 324.78] [color={rgb, 255:red, 0; green, 0; blue, 0 }  ][line width=0.75]    (10.93,-3.29) .. controls (6.95,-1.4) and (3.31,-0.3) .. (0,0) .. controls (3.31,0.3) and (6.95,1.4) .. (10.93,3.29)   ;

\draw (86,43) node [anchor=north west][inner sep=0.75pt]    {$P_{Spin}( X,g)$};
\draw (388,43) node [anchor=north west][inner sep=0.75pt]    {$P_{SO}( X,g)$};
\draw (257,151) node [anchor=north west][inner sep=0.75pt]    {$( X,g)$};
\draw (275,30) node [anchor=north west][inner sep=0.75pt]    {$\beta $};

\end{tikzpicture}

\end{center}
The spinor bundle $\Sigma_{g} X$ is then defined by 
$$\Sigma_{g} X:=P_{Spin}(X,g)\times_{\beta} \mathbb{S}_{n+1}.$$
where here we denoted by $\mathbb{S}_{n+1}$ the unique (up to isomorphism) irreducible complex $Cl_{n+1}$-module such that $Cl_{n+1}\otimes \C \equiv End_{\C}(\mathbb{S}_{n+1})$ as a $\C$-a]gebra, where $Cl_{n+1}$ denotes the Clifford algebra of $\R^{n+1}$. The Spinor bundle $\Sigma_{g} X$ carries a natural Hermitian structure and a metric connection induced by the Levi-Civita connection on $TX$, that will be denoted by $\nabla^{g}$. Moreover, there is a natural Clifford multiplication defined by the action of $TX$ on $\Sigma_{g} X$. Below, we summarize the main properties of the spinor bundle and we refer the reader to \cite{Gin,LM} for more background on spinorial geometry:
\begin{itemize}
\item For all $U, V \in C^{\infty}(TX)$ and $\psi \in C^{\infty}(\Sigma_{g} X)$ we have $U\cdot V \cdot \psi +V\cdot U \cdot \psi=-2g(U,V)\psi$. Here, "$\cdot$ " denotes the Clifford multiplication.
\item If $(\cdot,\cdot)$ denotes the Hermitian metric on $\Sigma_{g} X$, then for all $U \in C^{\infty}(TX)$ and $\psi, \phi \in C^{\infty}(\Sigma_{g} X)$ we have $(U\cdot \psi , \phi)=-(\psi, U\cdot \phi)$.
\item For all $\psi, \phi \in C^{\infty}(\Sigma_{g} X)$ and $U\in C^{\infty}(TX)$, we have $U(\psi,\phi)=(\nabla^{g}_{U}\psi,\phi)+(\psi,\nabla^{g}_{U} \phi)$.

\item For all $U, V \in C^{\infty}(TX)$ and $\psi \in C^{\infty}(\Sigma_{g} X)$ we have $\nabla_{U} (V\cdot \psi)=(\nabla_{U}V)\cdot \psi+V\cdot \nabla_{U}\psi$.
\end{itemize}
For the rest of the manuscript, we let $\langle \cdot, \cdot \rangle:= Re (\cdot,\cdot)$. Then $\langle \cdot, \cdot \rangle$ defines a metric on $\Sigma_{g} X$.\\
If the manifold $X$ has non-empty boundary $\partial X =M$, then $(M,h)$ inherits a natural Spin structure from $(X,g)$, where $h:=g|_{M}$. Indeed, we have the following properties:
\begin{itemize}
\item There exists a unitary isomorphism 
$$\mathcal{L}: \Sigma_{g} X |_{M}\to \left\{ \begin{array}{ll} 
\Sigma_{g} M \text{ if } n=2k\\
\Sigma_{h} M \oplus \Sigma_{h} M \text{ if } n=2k+1
\end{array}
\right.$$
\item If $\nu$ is the normal vector generated by the embedding $M\subset X$, then using the previous identification we have for all $U\in C^{\infty}(TM)$ and $\varphi \in C^{\infty}(\Sigma_{g}X)$,
 $$U\cdot \nu \cdot \varphi=\left\{\begin{array}{ll}
U\ \overline{\cdot}\ \varphi \text{ if }n=2k\\
\Big(U\ \overline{\cdot}\ \oplus -U\ \overline{\cdot}\ \Big)\varphi \text{ if } n=2k+1
\end{array}
\right.$$
where here $\overline{\cdot}$ denotes the Clifford multiplication on $\Sigma_{h}M$.
\item For all $U\in C^{\infty}(TM)$ and $\varphi \in C^{\infty}(\Sigma_{g}X)$
$$\nabla^{g}_{X}\varphi=\nabla^{h}_{X}\varphi+\frac{A(X)}{2}\cdot \nu \cdot \varphi,$$
where $A$ is the Weingarten endormorphism.
\end{itemize}

In particular, near the boundary $\partial X=M$, if we denote by $D_{g}$ the Dirac operator on $(X,g, \Sigma_{g}X)$, then we have
\begin{equation}\label{mean}
D_{g}=\nu\cdot \nabla^{\overline{g}}_{\nu}+\nu\cdot \Big( D_{2}-\frac{n}{2}H\Big),
\end{equation}
where $$D_{2}=\left\{\begin{array}{ll}
D_{h} \text{ if } n=2k\\
D_{h}\oplus -D_{h} \text{ if } n=2k+1.
\end{array}
\right.
$$
In what follows we will be using $D_{h}$ instead of $D_{2}$. It is important to keep in mind the properties of the action of the normal vector field $\nu$ via Clifford multiplication as follows: 

$$(\nu\cdot)^{2}=-1\quad \text{ and } \quad \nu \cdot D_{h}:=-D_{h}\nu\cdot.$$
It follows that the spinor bundle splits into $\Sigma_{h}M=\Sigma^{+}\oplus \Sigma^{-}$, with $\Sigma^{\pm}=\ker(\nu\cdot \pm i)$.
Now, we will put all the above in the setting of a Poincar\'{e}-Einstein manifold. So as in the introduction of Sub-Section 2.2, we consider a Poincar\'{e}-Einstein manifold  $(X^{n+1},g_{+})$ with $(M,h)$ being its boundary at infinity, that is, $(M,h)$ is the boundary of $(\overline{X}^{n+1},\overline{g}=x^{2}g_{+})$. This way, the spinor bundle $\Sigma_{h}M$ is well defined through the identification isomorphism $\mathcal{L}$. We will abuse the notation and identify $\Sigma_{\overline{g}}\overline{X}|_{M}$ with $\Sigma_{h}M$. The spectrum and resolvent of the Dirac operator on a Poincar\'{e}-Einstein manifold is well understood. We refer the reader to \cite{GMP,GMP2} for the full study and \cite{M,MM} for the general approach, but we give below the necessary properties that we will be using in our investigation:
\begin{proposition}[\cite{GMP,GMP2}]
The spectrum of the Dirac operator on a Poicare-Einstein Manifold $(X^{n+1},g_{+})$ is $\R$ and is absolutely continuous. The resolvent $R(\lambda):=(D^{2}_{g_{+}}+\lambda^{2})^{-1}$ extends meromorphically to $\C\setminus (-\N/2)$ with poles of finite multiplicity, except at $\lambda=0$. Moreover, $R(\lambda)$ maps $x^{\infty}C^{\infty}(\Sigma_{\overline{g}}\overline{X})$ to $x^{\frac{n}{2}+\lambda}C^{\infty}(\Sigma_{\overline{g}}\overline{X})$.
\end{proposition}

We will state below the central theorem that allows us to construct the scattering operator, along with its properties:
\begin{theorem}[\cite{GMP,GMP2}]\label{thmscat}
Let $\psi \in C^{\infty}(\Sigma_{h}M)$ and $\lambda \in \C\setminus -(\frac{1}{2}+\N)$ not a pole of $R(\lambda)$, then there exists  $\sigma \in C^{\infty}(\Sigma_{g_{+}}X)$ solution to 
\begin{equation}\label{eq1}
(D^{2}_{g_{+}}+\lambda^{2})\sigma =0.
\end{equation}
The solutions $\sigma$ is unique for $\Re(\lambda)\geq 0$. Moreover, there exists $\sigma_{\pm}\in C^{\infty}(\Sigma_{\overline{g}}\overline X^{n+1})$ such that $\sigma=x^{\frac{n}{2}-\lambda}\sigma_{+}+x^{\frac{n}{2}+\lambda}\sigma_{-}$ with $\sigma_{+}|_{x=0}=\psi$.
The two spinors $\sigma_{\pm}$ depend meromorphically on $\lambda \in \C\setminus -(\frac{1}{2}+\N)$.
\end{theorem}

The scattering operator $S(\lambda)$ is defined for $\lambda \in \C \setminus \frac{1}{2}+\N$ and $\Re(\lambda)\geq 0$ by $$S(\lambda) :\sigma_{+}(\lambda)|_{x=0}\mapsto \sigma_{-}(\lambda)|_{x=0}.$$
The main properties of $S$ can be summarized in the proposition below: 
\begin{proposition}[\cite{GMP}]
Under the assumptions of Theorem \ref{thmscat}, the scattering operator $S$ satisfies the following properties: 
\begin{itemize}
\item[i)] $S(\cdot)$ defines a holomorphic family of pseudo-differential operators on $\Omega=\{ \lambda \in \C \setminus \frac{1}{2}+\N; \Re(\lambda)\geq 0\}$ of order $2\lambda$. Moreover, it has simple poles at $\lambda= \frac{1}{2}+k$ for $k\in \N$ and there exists a constant $c_{k}$ and a differential operator $\D_{h}^{2k+1}$ such that
$$Res_{\lambda=\frac{1}{2}+k}S(\lambda)=c_{k}\D_{h}^{2k+1}.$$
\item[ii)] Under the conformal change $\hat{h}=e^{2w}h$ with $w \in C^{\infty}(M)$ and if $\hat{S}$ is the corresponding scattering operator for the metric $\hat{h}$, we have
$$\hat{S}(\lambda)=e^{-(\frac{n}{2}+\lambda)w}S(\lambda) e^{(\frac{n}{2}-\lambda)w}.$$
\item[iii)] For $\Re(\lambda)>0$, $S(\lambda)$ is essentially self-adjoint and has a principal symbol
$$\sigma_{pr}(S(\lambda))(\xi)=i 2^{-2\lambda}\frac{\Gamma(\frac{1}{2}-\lambda)}{\Gamma(\frac{1}{2}+\lambda)}|\xi|_{h}^{2\lambda-1}\nu\cdot \xi \cdot.$$
\end{itemize}
\end{proposition}

These properties are relatively similar to the ones for the scattering operator in the scalar case as in \cite{GZ}. The differential operators $\D_{h}^{2k+1}$ are the spinorial version of the GJMS operators as they appear in the residue part of the scattering operator. Notice that in contrast with the classical GJMS operators, the $\D_{h}^{2k+1}$ are odd dimensional with leading term $\nu\cdot D_{h}^{2k+1}$. We can then normalize them in a way that the leading term is $D_{h}^{2k+1}$ and that is by considering the operators $\overline{\D}^{2k+1}_{h}:=-\nu\cdot \D_{h}^{2k+1}$.
In fact, one has explicit expressions for $\overline{\D}_{h}^{2k+1}$ for $k=0,1,2$ ( see \cite{GMP2,F,F2}) and for all $k$ if $(M,h)$ is Einstein (see \cite{FKS}). Indeed, $\overline{\D}_{h}^{1}=D_{h}$ and
$$\overline{\D}_{h}^{3}=D_{h}^{3}\psi-\frac{2 Ric^{h}_{ij}e_{i}\cdot \nabla_{j}^{h}\psi}{n-2}-\frac{R_{h}}{(n-1)(n-2)}D_{h}\psi-dR_{h}\cdot\psi.$$
If $(M,h)$ is Einstein, then one has an explicit expression for $\overline{\D}_{h}^{2k+1}$ as shown in \cite{FKS}:
$$\overline{\D}_{h}^{2k+1}=D_{h}\prod_{j=1}^{k}\Big(D_{h}^{2}-j^{2}\frac{R_{h}}{n(n-1)}\Big).$$
One also can express the operators $\D_{h}^{2k+1}$ using the equation (\ref{eq1}), similar to \cite{GZ}:
\begin{theorem}[\cite{GMP2,F}]
For $\lambda =\frac{1}{2}+k$ with $k\in \N$ ($k\leq \frac{n}{2}$ if $n$ is even) and $\psi \in C^{\infty}(\Sigma_{h}M)$, there exists $\sigma \in C^{\infty}(X^{n+1},g_{+})$ solution of
$$(D^{2}_{g_{+}}+\lambda^{2})\sigma =0.$$ Near the boundary $\partial \overline{X}=M$, there exists $F, G\in C^{\infty}(\Sigma_{\overline{g}}\overline{X})$ such that $$\sigma=x^{\frac{n}{2}-\lambda}F+x^{\frac{n}{2}+\lambda}\ln(x)G,$$
with $F|_{x=0}=\psi$ and $G|_{x=0}=Res_{\lambda=\frac{1}{2}+k}S(\lambda)\psi=c_{k}\D_{h}^{2k+1}\psi$. The spinor $F$ is unique up to order $O(x^{2k+1})$ and $G$ is unique up to order $O(x)$.
\end{theorem}
We will give here some details about the formal construction of the solution $\sigma$ of (\ref{eq1}) in order to compare with the scalar case. Indeed, if we define 
$$\mathcal{A}:=\{\varphi \in C^{\infty}(\Sigma_{\overline{g}}\overline{X}^{n+1}); \varphi=\sum_{j=0}^{\infty}x^{j}\varphi_{j}; \nabla^{\overline{g}}_{\nu}\varphi_{j}=0\},$$
and 
$$\mathcal{A}^{\pm}:=\{\varphi \in \mathcal{A}; \varphi_{2j}\in C^{\infty}(\Sigma^{\pm}M), \varphi_{2j+1}\in C^{\infty}(\Sigma^{\mp}M)\}.$$
The solution $\sigma$ of $(D^{2}_{g_{+}}+\lambda^{2})\sigma=0$ can be then decomposed as $\sigma=\sigma^{+}+\sigma^{-}$ where $\sigma^{\pm}$ is a solution of the equation \begin{equation}\label{ord1}
(D_{g_{+}}\pm i\lambda)\sigma^{\pm}=0,
\end{equation} and near the boundary $\sigma^{\pm}=x^{\frac{n}{2}-\lambda}\sigma_{+}^{\pm}+x^{\frac{n}{2}+\lambda}\sigma_{-}^{\pm}$ such that $\sigma_{+}^{\pm}\in \mathcal{A}^{\pm}$ and $\sigma_{-}^{\pm}\in \mathcal{A}^{\mp}$. In particular,
$$S(\lambda)=S^{+}(\lambda)+S^{-}(\lambda), \text{ and } S^{\pm}(\lambda):C^{\infty}(\Sigma^{\pm}M)\to C^{\infty}(\Sigma^{\mp}M).$$
Notice that, in the scalar case (Theorem \ref{thmscal}), the expansion of $F$ and $G$ is only in even powers of $x$, but for the spinorial case, the expansion is in odd and even powers of $x$.\\

To finish this subsection, we introduce the wighted Dirac and Laplace operators. These operators will be relevant for our study later on. For more details regarding the study of spinors on weighted manifolds, we refer the reader to \cite{BO}. Given a smooth function $f\in C^{\infty}(X^{n+1},\overline{g})$ we define the weighted Dirac and Laplacian as follows:
$$D_{\overline{g}}^{f}=D_{\overline{g}}-\frac{1}{2}\nabla^{\overline{g}}f \cdot$$
and
$$-\Delta_{\overline{g}}^{f}=-\Delta_{\overline{g}}+\nabla^{\overline{g}}_{\nabla^{\overline{g}}f}.$$
The main properties of these operators can be summarized below:
\begin{proposition}\label{int}
Given $\psi, \varphi \in C^{\infty}(\Sigma_{\overline{g}}\overline{X}^{n+1})$, then
$$\int_{X}\langle D_{\overline{g}}^{f}\psi ,\varphi \rangle e^{-f}\ dv_{\overline{g}}=\int_{X}\langle \psi, D_{\overline{g}}^{f}\varphi\rangle e^{-f}\ dv_{\overline{g}}+\int_{M} \langle \nu \cdot \psi,\varphi \rangle e^{-f}\ dv_{h},$$
and
$$\int_{X}\langle -\Delta_{\overline{g}}^{f}\psi,\varphi \rangle e^{-f}\ dv_{\overline{g}}=\int_{X}\langle \nabla^{\overline{g}}\psi,\nabla ^{\overline{g}}\varphi\rangle e^{-f}\ dv_{\overline{g}}-\int_{M}\langle \nabla^{\overline{g}}_{\nu} \psi,\varphi \rangle e^{-f} \ dv_{h}.$$
Moreover,
$$D_{\overline{g}}^{2}=-\Delta_{\overline{g}}+\frac{R^{f}_{\overline{g}}}{4},$$
where $R_{\overline{g}}^{f}=R_{\overline{g}}+2\Delta_{\overline{g}}f-|\nabla^{\overline{g}}f|^{2}$.
\end{proposition}
In fact, the operator $D_{\overline{g}}^{f}$ is the Dirac operator for the modified connection $\nabla^{\overline{g},f}$ defined for $U\in C^{\infty}(TX)$ and $\psi \in C^{\infty}(\Sigma_{\overline{g}} X)$ by
$$\nabla^{\overline{g},f}_{U}\psi=\nabla_{U}^{\overline{g}}\psi -\frac{1}{2}(\nabla^{\overline{g}}_{U}f)\psi.$$
One can also express $D_{\overline{g}}^{f}$ as 
$$D_{\overline{g}}^{f}=e^{\frac{f}{2}}D_{\overline{g}}e^{-\frac{f}{2}}.$$

\subsection{The Conformal Fractional Dirac Operator and The Spinorial Yamabe Problem}

\begin{definition}
Assume that $\lambda>0$ ($\lambda \leq \frac{n}{2}$ if $n$ is even) with $\lambda\not \in \frac{1}{2}+\N$ and let $d_{\lambda}=2^{2\lambda}\frac{\Gamma(\frac{1}{2}+\lambda)}{\Gamma(\frac{1}{2}-\lambda)}$. The conformal fractional Dirac operator is defined by $$\D^{2\lambda}_{h}:=d_{\lambda}S(\lambda).$$
The geometric conformal fractional Dirac operator is then defined by
$$\overline{\D}_{h}^{2\lambda}=-\nu\cdot \D_{h}^{2\lambda}.$$
\end{definition}
The fractional Dirac operator $\D_{h}^{2\lambda}$ and its geometric counterpart $\overline{\D}_{h}^{2\lambda}$ are essentially self-adjoint and conformally invariant. That is, if $\hat{h}=e^{2u}h$ then
$$\D_{\hat{h}}^{2\lambda}=e^{-(\frac{n}{2}+\lambda)u}\D_{h}^{2\lambda}e^{(\frac{n}{2}+\lambda)u}.$$
This last property follows directly from the conformal invariance of the scattering operator $S(\lambda)$. For the geometric operator, the conformal invariance is deduced from the fact that $\nu\cdot$ is conformally invariant at the boundary of $\overline{X}^{n+1}$. To prove the self-adjointness of $\D_{h}^{2\lambda}$, we consider two solutions $\sigma_{1}$ and $\sigma_{2}$ of (\ref{eq1}) with $x^{\lambda-\frac{n}{2}}\sigma_{i}|_{x=0}=\psi_{i}\in C^{\infty}(\Sigma_{h}M)$, $i=1,2$. Then we have

$$F.P. \int_{x>\varepsilon}\langle \nabla \sigma_{1},\nabla \sigma_{2}\rangle +\Big(\lambda^{2}-\frac{n(n+1)}{2}\Big)\langle \sigma_{1},\sigma_{2}\rangle \ dv_{g_{+}}=\frac{n}{d_{\lambda}}\int_{M} \langle \D_{h}^{2\lambda}\psi_{1},\psi_{2}\rangle \ dv_{h},$$
where $F.P.$ denotes the finite part of the integral when $\varepsilon\to 0$. Now the symmetry of the left hand-side implies the symmetry of the right-hand-side. 
Moreover, from the expression of the principal symbol of $S(\lambda)$, we see that $\D_{h}^{2\lambda}=|D_{h}|^{2\lambda-1}\nu\cdot D_{h}+l.o.t.$\\

Notice that if we set $\tilde{d}_{\lambda}=\partial_{\lambda}d_{\lambda}$ then we have
$$\lim_{\lambda\to \frac{1}{2}+k}\D_{h}^{2\lambda}=\tilde{d}_{\frac{1}{2}+k}Res_{\lambda=\frac{1}{2}+k}S(\lambda)=\D_{h}^{2k+1}.$$
We provide below two explicit expressions of $\D_{h}^{2\lambda}$ that are proved in the Appendix:
\begin{proposition}\label{propex}
In the case of $\R^{n}$ with its flat metric and standard spin structure, we have
$$\D^{2\lambda}_{\R^{n}}=|D_{\R^{n}}|^{2\lambda-1}\nu\cdot D_{\R^{n}}, \text{ for $\lambda\geq 0$}.$$
On the other hand, for the standard sphere $(S^{n},h_{0})$ the conformal fractional Dirac operator takes the form:
$$\D_{h_{0}}^{2\lambda}=\frac{\Gamma\Big(|\nu\cdot D_{h_{0}}|+\frac{1}{2}+\lambda\Big)}{\Gamma\Big(|\nu\cdot D_{h_{0}}|+\frac{1}{2}-\lambda\Big)}\frac{\nu\cdot D_{h_{0}}}{|\nu \cdot D_{h_{0}}|}, \text{ for all $\lambda>0$}.$$
Moreover,
$$\D^{2(\lambda+1)}_{h_{0}}=\Big(D_{h_{0}}^{2}-(\lambda+\frac{1}{2})^{2}\Big)\D_{h_{0}}^{2\lambda}.$$
\end{proposition}

Using the properties above one can state a fractional version of the Spinorial Yamabe Problem. Indeed, we define the fractional  Spinorial Yamabe constant $\lambda_{2\lambda}^{+}(M,[h])$ by
$$\lambda_{2\lambda}^{+}(M,[h]):=\inf_{\hat{h}\in [h]}|\lambda_{1}(\overline{\D}_{\hat{h}}^{2\lambda})|Vol(\hat{h})^{\frac{2\lambda}{n}}.$$
We claim the following:
\begin{theorem}\label{thmyam}
\begin{equation}\label{eqyam}
\lambda_{2\lambda}^{+}(M,[h])=\inf_{\varphi \in L^{\frac{2n}{n+2\lambda}}(\Sigma_{h}M)\setminus \{0\}} \frac{\|\varphi\|_{L^{\frac{2n}{n+2\lambda}}}^{2}}{|\int_{M}\langle |\overline{\D}_{h}^{2\lambda}|^{-1}\varphi,\varphi \rangle \ dv_{h}|},
\end{equation}
and if the inf is achieved then there exists $\psi$ that satisfies the equation \begin{equation}\label{eqyam1}
\overline{\D}_{h}^{2\lambda}\psi=\lambda^{+}_{2\lambda}|\psi|^{\frac{4\lambda}{n-2\lambda}}\psi.
\end{equation}
Moreover, there exists $c=c([h])>0$ such that $$\lambda_{2\lambda}^{+}(M,[h])\geq c.$$
\end{theorem}
{\it Proof:}
The proof of these results is the same as in the case of the classical Dirac operator proved in \cite{Amm2}. The proof of the conformal lower bound follows the idea of Lott \cite{L} that was then generalized to the case when $D_{h}$ is not invertible in \cite{Amm2}. For simplicity, we will assume that $\overline{\D}_{h}^{2\lambda}$ is invertible. The generalization is clear when $\ker \overline{\D}_{h}^{2\lambda}\not=\{0\}$.\\
Assume that $\hat{h}=f^{2}h_{0}$, and let $0<\lambda_{1}(\hat{h})$ be the first eigenvalue of the operator $\Big(\overline{\D}_{\hat{h}}^{2\lambda}\Big)^{2}$. For $\psi \in L^{2}(\Sigma_{\hat{h}}M)$ we let $\varphi=f^{\frac{(n+2\lambda)}{2}}\psi$. Then we have
\begin{align}
\lambda_{1}^{-\frac{1}{2}}(\hat{h})&=\sup_{\psi \in L^{2}(\Sigma_{\hat{h}}M)\setminus \{0\}}\frac{|\int_{M}\langle \Big(\overline{\D}_{\hat{h}}^{2\lambda}\Big)^{-1}\psi,\psi \rangle \ dv_{\hat{h}}|}{\int_{M}|\psi|^{2}\ dv_{\hat{h}}}\label{sup}\\
&=\sup_{\varphi \in L^{2}(\Sigma_{h}M)\setminus \{0\}}\frac{|\int_{M}\langle \Big(\overline{\D}_{\hat{h}}^{2\lambda}\Big)^{-1}f^{-\frac{(n+2\lambda)}{2}}\varphi,f^{-\frac{(n+2\lambda)}{2}}\varphi \rangle f^{n}\ dv_{h}|}{\int_{M}|\varphi|^{2}f^{-2\lambda}\ dv_{h}}\notag\\
&=\sup_{\varphi \in L^{2}(\Sigma_{h}M)\setminus \{0\}}\frac{|\int_{M}\langle \Big(\overline{\D}_{h}^{2\lambda}\Big)^{-1}\varphi,\varphi \rangle \ dv_{h}|}{\int_{M}|\varphi|^{2}f^{-2\lambda}\ dv_{h}}.\notag
\end{align}
On the other hand, from H\"{o}lder's inequalities, we have
$$\int_{M}|\varphi|^{\frac{2n}{n+2\lambda}}\ dv_{h}\leq \Big(\int_{M}|\varphi|^{2}f^{-2\lambda}\ dv_{h}\Big)^{\frac{n}{n+2\lambda}}Vol(\hat{h})^{\frac{2\lambda}{n+2\lambda}},$$
with equality if and only if $f=c|\varphi|^{\frac{2}{n+2\lambda}}$.
Keeping in mind that the sup in (\ref{sup}), can be taken for $\varphi \in L^{\frac{2n}{n+2\lambda}}(\Sigma_{h}M)$, we have
\begin{align}
\lambda_{2\lambda}^{+}(M,[h])&=\inf_{\{f\in C^{\infty}(M),f>0\}}\inf_{\{\varphi \in L^{\frac{2n}{n+2\lambda}}(\Sigma_{h}M)\setminus \{0\}\}}\frac{Vol(\hat{h})^{\frac{2\lambda}{n}}\int_{M}f^{-2\lambda}|\varphi|^{2}\ dv_{h}}{|\int_{M} \langle |\overline{\D}_{h}^{2\lambda}|^{-1}\varphi,\varphi \rangle \ dv_{h}|}\notag\\
&\geq \inf_{\varphi \in L^{\frac{2n}{n+2\lambda}}(\Sigma_{h}M)\setminus \{0\}} \frac{\|\varphi\|_{L^{\frac{2n}{n+2\lambda}}}^{2}}{|\int_{M}\langle |\overline{\D}_{h}^{2\lambda}|^{-1}\varphi,\varphi \rangle \ dv_{h}|}.\notag
\end{align}
The equality then follows by taking $f=|\varphi|^{\frac{2}{n+2\lambda}}$ if $|\varphi|\not=0$ or using an approximation $f_{k}$ so that $$\|f_{k}-\max\{\frac{1}{k},|\varphi|^{\frac{2}{n+2\lambda}}\}\|_{C^{0}}\leq e^{-k}, \text{ for $k\in \N$ }.$$
For a minimizer $\varphi$ of (\ref{eqyam}), we can easily deduce by setting $\psi=(\overline{\D}^{2\lambda}_{h})^{-1}\varphi$ that 
$$\overline{\D}_{h}^{2\lambda}\psi=\lambda_{2\lambda}^{+}(M,[h])|\psi|^{\frac{4\lambda}{n-2\lambda}}\psi.$$
Consider now the operator $\mathcal{K}=(-\Delta_{h}+1)^{\frac{\lambda}{2}}|\overline{\D}_{h}^{2\lambda}|^{-\frac{1}{2}}$. Then $\mathcal{K}$ is a pseudo-differential operator of order zero. In particular, it maps $L^{p}$ to $L^{p}$ for all $p\in (1,\infty)$. Hence, $\mathcal{K}$ is a well defined bounded operator from $L^{\frac{2n}{n+2\lambda}}(\Sigma_{h}M)\to L^{\frac{2n}{n+2\lambda}}(\Sigma_{h}M)$. Therefore, by elliptic regularity of $(-\Delta_{h} +1)^{\frac{\lambda}{2}}$, we have that
$|\overline{\D}_{h}^{2\lambda}|^{-\frac{1}{2}}$ is  a bounded operator from $L^{\frac{2n}{n+2\lambda}}(\Sigma_{h}M)$ to $W^{\lambda,\frac{2n}{n+2\lambda}}(\Sigma_{h}M)$. Finally, using the Sobolev embedding $W^{\lambda,\frac{2n}{n+2\lambda}}(\Sigma_{h}M)\hookrightarrow L^{2}(\Sigma_{h}M)$, we see that there exists $a>0$ such that
$$|\int_{M}\langle |\overline{\D}_{h}^{2\lambda}|^{-1}\varphi,\varphi \rangle \ dv_{h}|\leq a \|\varphi\|_{L^{\frac{2n}{n+2\lambda}}}^{2}, \text{ for all $\varphi \in L^{\frac{2n}{n+2\lambda}}(\Sigma_{h}M)$}.$$
Therefore,
$$\lambda_{2\lambda}^{+}(M,[h])\geq \frac{1}{a}>0.$$
\hfill $\Box$

\begin{example}\label{ex1}
From the expression of $\D_{h_{0}}^{2\lambda}$ on the sphere $(S^{n},h_{0})$ we see from Appendix A that 
$$\lambda_{1}(\D^{2\lambda}_{h_{0}})=\frac{\Gamma(\frac{n}{2}+\frac{1}{2}+\lambda)}{\Gamma(\frac{n}{2}+\frac{1}{2}-\lambda)},$$
and it corresponds to the eigenspinor $\hat{\varphi}_{1,+}=(1-\nu\cdot)\varphi_{-1}$, where $\varphi_{-1}$ is a $\frac{1}{2}$-Killing spinor. We also point out that the operator $\overline{\D}^{2\lambda}_{h_{0}}$ has the same spectrum and the same first eigenvalue $\lambda_{1}(\overline{\D}_{h_{0}}^{2\lambda})$ with eigenspinor $\varphi_{1}$. In particular, since $|\varphi_{1}|$ is constant, we can normalize it so that $|\varphi_{1}|=1$. Hence, $\varphi_{1}$ satisfies equation (\ref{eqyam1}), that is,
$$\overline{\D}_{h_{0}}^{2\lambda}\varphi_{1}=\lambda_{1}(\overline{\D}_{h_{0}}^{2\lambda})|\varphi_{1}|^{\frac{4\lambda}{n-2\lambda}}\varphi_{1}.$$

Consider now, the stereographic projection $\pi:S^{n}\setminus \{N\}\to \R^{n}$, where $N\in S^{n}$ is the north pole. Then we have $(\pi^{-1})^{*}h_{0}=f^{2}h_{\R^{n}}$, where $f(x)=\frac{2}{1+|x|^{2}}$. On the other hand, we recall that the Killing spinors on the sphere take the form
$$\varphi_{\pm 1}=\pi^{*}(f^{\frac{1}{2}}(1\mp x\cdot)\Phi_{0}),$$
where $\Phi_{0}$ is a constant spinor on $\R^{n}$ and, here, $\pi^{*}$ is the identification map between conformally equivalent spinor bundles introduced in (\ref{iden}). Therefore, if we set $\psi= f^{\frac{n+1-2\lambda}{2}}(1-x\cdot) \Phi_{0}$, we have 
\begin{equation}\label{test}
\psi=f^{\frac{n-2\lambda}{2}}(\pi^{-1})^{*}\varphi_{1}.
\end{equation}
Therefore, we have
$$\overline{\D}_{\R^{n}}^{2\lambda}\psi=\lambda_{1}(\overline{\D}^{2\lambda}_{h_{0}})|\psi|^{\frac{4\lambda}{n-2\lambda}}\psi.$$
Indeed, from the conformal invariance property of $\overline{\D}^{2\lambda}$ we have
\begin{align}
\overline{\D}_{\R^{n}}^{2\lambda}(\psi)&=f^{\frac{n+2\lambda}{2}}(\pi^{-1})^{*}\Big(\overline{\D}^{2\lambda}_{h_{0}}\pi^{*}(f^{\frac{1}{2}}(1-x\cdot) \Phi_{0})\Big)\notag\\
&=f^{\frac{n+2\lambda}{2}}(\pi^{-1})^{*}\overline{\D}^{2\lambda}_{h_{0}}(\varphi_{1})\notag\\
&=\lambda_{1}(\overline{\D}^{2\lambda}_{h_{0}})f^{\frac{n+2\lambda}{2}}(\pi^{-1})^{*}(\varphi_{1})\notag\\
&=\lambda_{1}(\overline{\D}^{2\lambda}_{h_{0}})f^{2\lambda}\psi.\notag
\end{align}
But $|\psi|=f^{\frac{n-2\lambda}{2}}$. Hence,
$$\overline{\D}_{\R^{n}}\psi=\lambda_{1}(\overline{\D}^{2\lambda}_{h_{0}})|\psi|^{\frac{4\lambda}{n-2\lambda}}\psi.$$
\end{example}

\begin{remark}
The results in Theorem \ref{thmyam} leave many of open questions. The most natural ones are
\begin{itemize}
\item[i)] Can we find an explicit lower bound for the constant $c$? Notice that for the case of the Dirac operator ($\lambda=\frac{1}{2}$) one has the Hijazi's inequality, which provides a lower bound on $\lambda^{+}_{1}$ in terms of the Yamabe constant. Threfore, can one find a lower bound that might involve the fractional Yamabe constant or the fractional $Q$-curvature?
\item[ii)] Can one prove an Aubin-type inequality for $\lambda_{2\lambda}^{+}$?  That is $\lambda_{2\lambda}^{+}(M,[h])\leq \lambda_{2\lambda}^{+}(S^{n},[h_{0}])$?
\end{itemize}
The second question is traditionally proved by using a "good" test spinor in the functional $I$ (or another more adapted functional). The natural choice for the test spinor would be a graft of the "standard bubble", which is the solution of (\ref{eqyam1}) in $\R^{n}$. That is, as we have seen in Example \ref{ex1}, the candidate for the test spinor is $$\psi=\Big(\frac{2}{1+|x|^{2}}\Big)^{\frac{n-2\lambda+1}{2}}(1-x\cdot)\Phi_{0}.$$
 But, as we see from Proposition \ref{propex}, the expression of the fractional Dirac operator is not computationally friendly and we need to find a better way that would allow us to perform computations using $\D_{h}^{2\lambda}$. This leads us to the next part of the paper, which is transforming the problem into a Dirichlet-to-Neuman type problem, or what is known in the literature by the Caffarelli-Silvestre extentions.
\end{remark}
\section{The Extension Problem}
In this section, we focus on the extension problem corresponding to the conformal compactification of $(X^{n+1},g_{+})$. As we will see, different compactifications can lead to different kinds of precision. 
\begin{theorem}\label{thmext1}
Let $\lambda \in (0,\frac{n}{2}]$ and $\lambda \not \in \frac{1}{2}+\N$. Fix $\psi \in C^{\infty}(\Sigma_{h}M)$ and $f=\ln(x^{2\lambda-1})$, then $\Psi$ satisfies
\begin{equation}\label{BVP}
\left\{
\begin{array}{ll}
-\Delta^{f}_{\overline{g}}\Psi+\frac{|\nabla^{\overline{g}} x|}{x}\Big[\nabla^{\overline{g}}_{\nu}+\nu\cdot D_{\overline{g}}\Big]\Psi +E_{\lambda}(x)\Psi=0\\
\Psi|_{x=0}=\psi,
\end{array}
\right.
\end{equation}
with $E_{\lambda}(x)=\Big[\lambda\frac{\Delta_{\overline{g}}x}{x}+\lambda^{2}\Big(\frac{1-|\nabla^{\overline{g}} x|^{2}}{x^{2}}\Big)+\frac{R_{\overline{g}}}{4}\Big]$, if and only if $\sigma=x^{\frac{n}{2}-\lambda}\Psi$ satisfies 
$$\left\{\begin{array}{ll}
(D^{2}_{g_{+}}+\lambda^{2})\sigma=0\\
\sigma=x^{\frac{n}{2}-\lambda}\sigma_{+}+x^{\frac{n}{2}+\lambda}\sigma_{-}.
\end{array}
\right.
$$
with $\sigma_{+}|_{x=0}=\psi$. In this case we have 
\begin{itemize}
\item If $0<\lambda<\frac{1}{2}$, then
$$\D_{h}^{2\lambda}\psi=\frac{d_{\lambda}}{2\lambda}\lim_{x\to 0}x^{1-2\lambda}\partial_{x}\Psi.$$
\item If $\frac{1}{2}<\lambda<1$, then
$$\D_{h}^{2\lambda}\psi=\frac{d_{\lambda}}{2\lambda}\lim_{x\to 0} x^{1-2\lambda} (\partial_{x}\Psi-\frac{1}{1-2\lambda}\nu\cdot D_{h}\psi).$$
\end{itemize}
\end{theorem}
{\it Proof:}
We first notice that from the conformal invariance of $D_{g_{+}}$ on $\overline{X}^{n+1}$ we have
$$D_{g_{+}}\sigma=x^{\frac{n}{2}+1}D_{\overline{g}}x^{-\frac{n}{2}}\sigma.$$
In particular for $\sigma =x^{\frac{n}{2}-\lambda}\Psi$, we have
\begin{align}
(D_{g_{+}}^{2}+\lambda^{2})\sigma&=x^{\frac{n}{2}+1}D_{\overline{g}}(xD_{\overline{g}}x^{-\lambda}\Psi)+\lambda^{2}x^{\frac{n}{2}-\lambda}\Psi\notag\\
&=x^{\frac{n}{2}+1}D_{\overline{g}}(-\lambda x^{-\lambda}\nabla^{\overline{g}}x \cdot \Psi +x^{-\lambda+1}D_{\overline{g}}\Psi)+\lambda^{2}x^{\frac{n}{2}-\lambda}\Psi\notag\\
&=x^{\frac{n}{2}+1}\Big(-\lambda^{2}x^{-\lambda-1}|\nabla^{\overline{g}}x|^{2}\Psi-\lambda x^{-\lambda}D_{\overline{g}}(\nabla^{\overline{g}}x\cdot \Psi)+(-\lambda+1)x^{-\lambda}\nabla^{\overline{g}}x\cdot D_{\overline{g}}\Psi\notag\\
&\quad+x^{-\lambda+1}D^{2}_{\overline{g}}\Psi\Big)+\lambda^{2}x^{\frac{n}{2}-\lambda}\Psi.\notag
\end{align}
On the other hand,
\begin{align}
D_{\overline{g}}(\nabla^{\overline{g}}x\cdot \Psi)&=D_{\overline{g}}^{2}(x\Psi)-D_{\overline{g}}(xD_{\overline{g}})\notag\\
&=xD_{\overline{g}}^{2}\Psi -2\nabla_{\nabla x}\Psi-(\Delta_{\overline{g}}x)\Psi-\nabla x \cdot D_{\overline{g}}-xD_{\overline{g}}^{2}\Psi\notag\\
&=-2\nabla_{\nabla x}\Psi-(\Delta_{\overline{g}}x)\Psi-\nabla x \cdot D_{\overline{g}}.\notag
\end{align}
Therefore,
\begin{align}
(D_{g_{+}}^{2}+\lambda^{2})\sigma&=\Big[ [\lambda^{2}(1-|\nabla^{\overline{g}}x|^{2})+\lambda x\Delta_{\overline{g}}x]\Psi +x\nabla^{\overline{g}}x\cdot D_{\overline{g}}\Psi+2\lambda x\nabla^{\overline{g}}_{\nabla^{\overline{g}}x}\Psi+x^{2}D_{\overline{g}}^{2}\Psi\Big]\notag\\
&=x^{\frac{n}{2}-\lambda+2}\Big[ -\Delta_{f}\Psi+\frac{|\nabla^{\overline{g}}x|}{x}(\nabla^{\overline{g}}_{\nu}+\nu\cdot D_{\overline{g}})\Psi +\Big(\lambda\frac{\Delta_{\overline{g}}x}{x}+\lambda^{2} \Big(\frac{1-|\nabla^{\overline{g}}x|^{2}}{x^{2}}\Big)+\frac{R_{\overline{g}}}{4}\Big)\Psi\Big],\notag
\end{align}
which proves the first claim of the theorem. Now, assume that $\lambda \in (0,\frac{1}{2})$, then we have 
$$\Psi =x^{-\frac{n}{2}+\lambda}\sigma =\psi + x^{2\lambda}S(\lambda)\psi+o(x^{2\lambda}).$$
Hence, $$x^{1-2\lambda}\partial_{x}\Psi=2\lambda S(\lambda)\psi +o(1),$$
and this finishes the proof of the first point. If $\lambda \in (\frac{1}{2},1)$, then we have
$$\Psi =\psi +\frac{\nu\cdot D_{h}\psi}{1-2\lambda} x+x^{2\lambda}S(\lambda)\psi +o(x^{2\lambda}).$$
Thus $\partial_{x}\Psi -\frac{1}{1-2\lambda}\nu\cdot D_{h}\psi =2\lambda x^{2\lambda-1}S(\lambda)\psi$, which finishes the proof.
\hfill $\Box$\\
\begin{remark}
Notice that in the previous theorem, there is a main difference in comparison to the scalar case, and that is in the second term of (\ref{BVP}). Nevertheless one should observe that
$$\nabla_{\nu}^{\overline{g}}+\nu\cdot D_{\overline{g}}=D_{h_{x}}-\frac{n}{2}H_{x},$$
where $D_{h_{x}}$ is the Dirac operator on $(M,h_{x})$ and $H_{x}$ is the mean curvature of $(M,h_{x})$ as a boundary of $X^{n+1}$ at a fixed level set $x$. In particular, in the case of the flat space $\R^{n}$, one can consider the extension on $\R^{n+1}_{+}:=\{(x_{1},\cdots,x_{n},t)\in \R^{n+1};t>0\}$ to get
\begin{equation}\label{PBE0}
\left\{\begin{array}{ll}
div_{\R^{n+1}}(t^{1-2\lambda}\nabla\Psi)+t^{-2\lambda}D_{\R^{n}}\Psi =0\\
\Psi|_{x=0}=\psi,
\end{array}
\right.
\end{equation}
and if $		\lambda<\frac{1}{2}$, we have
$$\D_{h}^{2\lambda}\psi=\frac{d_{\lambda}}{2\lambda}\lim_{t\to 0}t^{1-2\lambda}\partial_{t}\Psi.$$
\end{remark}

When $\lambda$ is large, one can still recover the operator $\D_{h}^{2\lambda}$ as follows:

\begin{corollary}
Let $k=[2\lambda]$ and $m_{0}=1+k-2\lambda$ and $\lambda \in (0,\frac{n}{2})$, $\lambda \not \in \frac{1}{2}+\N$. Assume that $\Psi$ solves (\ref{BVP}). Then we have 
$$\D_{h}^{2\lambda}\psi=\frac{d_{\lambda}}{2\lambda(2\lambda-1)\cdots (2\lambda-k)}\lim_{x\to 0}x^{m_{0}}\partial^{k+1}_{x}\Psi.$$
\end{corollary}
{\it Proof:}

This again follows from the expansion $\Psi =\psi+\sum_{i=1}^{k}x^{i}\psi_{i}+x^{2\lambda}S(\lambda)\psi +o(x^{2\lambda})$, where $k=[2\lambda]$. Thus
$$x^{1+k-2\lambda}\partial_{x}^{k+1}\Psi =2\lambda(2\lambda-1)\cdots (2\lambda-k) S(\lambda)\psi+o(1).$$
\hfill $\Box$

Relaxing the assumption on the smoothness of the defining function, one could use a $\lambda$-admissible defining function $\rho$, in the sense of Definition \ref{admis} instead. We will consider in particular, the defining function that takes the form $\rho=x(1+\frac{Q_{h}^{\lambda}}{c_{\lambda}}x^{2\lambda}+o(x^{2\lambda}))$. With this in mind, we can state the following:
\begin{theorem}\label{thmext2}
Let $\lambda \in (0,\frac{1}{2})$ and consider the $\lambda$-admissible defining function $\rho$ described above. Let $\tilde{g}=\rho^{2}g_{+}$ and assume that $\Psi$ solves
$$\left\{
\begin{array}{ll}
-\Delta_{\tilde{g}}^{f}\Psi+\frac{|\nabla^{\tilde{g}} \rho|}{\rho}\Big[\nabla^{\tilde{g}}_{\nu}+\nu\cdot D_{\tilde{g}}\Big]\Psi +\frac{1}{4(1-2\lambda)}R_{\tilde{g}}\Psi=0,\\
\Psi|_{\rho=0}=\psi.
\end{array}
\right.
$$
Then $$\lim_{\rho\to 0} \frac{d_{\lambda}}{2\lambda}\rho^{1-2\lambda}\partial_{\rho}\Psi=\D_{h}^{2\lambda}\psi-(\frac{n}{2}-\lambda)\frac{d_{\lambda}}{c_{\lambda}}Q_{h}^{\lambda} \psi.$$
\end{theorem}
{\it Proof:} 
Since $\rho$ is the standard $\lambda$-admissible defining function, we see from (\ref{curv}), that $E_{\lambda}(\rho)=\frac{1}{4(1-2\lambda)}R_{\tilde{g}}$. Moreover, from the proof of Theorem \ref{thmext1}, we have that $\sigma =\rho^{\frac{n}{2}-\lambda}\Psi$ is the solution of (\ref{eq1}). Thus 
\begin{align}
\Psi&=\Big(\frac{x}{\rho}\Big)^{\frac{n}{2}-\lambda}\Big(\psi +x^{2\lambda}S(\lambda)\psi +o(x^{2\lambda}\Big)\notag\\
&=\Big(1-(\frac{n}{2}-\lambda)\frac{Q_{h}^{\lambda}}{c_{\lambda}}x^{2\lambda}+o(x^{2\lambda})\Big)\Big(\psi +x^{2\lambda}S(\lambda)\psi +o(x^{2\lambda})\Big)\notag\\
&=\psi+\Big(S(\lambda)\psi -(\frac{n}{2}-\lambda)\frac{Q_{h}^{\lambda}}{c_{\lambda}}\psi \Big)x^{2\lambda}+o(x^{2\lambda})\notag\\
&=\psi +\Big(S(\lambda)\psi -(\frac{n}{2}-\lambda)\frac{Q_{h}^{\lambda}}{c_{\lambda}}\psi \Big)\rho^{2\lambda}+o(\rho^{2\lambda}).\notag
\end{align}
Thus,
$$\rho^{1-2\lambda}\partial_{\rho}\Psi=2\lambda \Big(S(\lambda)\psi -(\frac{n}{2}-\lambda)\frac{Q_{h}^{\lambda}}{c_{\lambda}}\psi\Big)+o(1),$$
which finishes the proof.
\hfill$\Box$

This last connection between $\D_{h}^{2\lambda}$ and $Q_{h}^{\lambda}$ can be seen as a fractional version of the link between the Dirac operator and the mean curvature at the boundary as in the identity (\ref{mean}).\\

The previous results can be formulated using a first order operator as it was the case for the definition of the scattering operator using equation (\ref{ord1}).
\begin{proposition}
Let $f_{1}=\ln(x^{2\lambda})$ and assume that $\lambda \not \in \frac{1}{2}+\N$. Let $\psi^{\pm}\in C^{\infty}(\Sigma^{\pm}M)$, then $\sigma^{\pm}$ is the unique solution of
$$(D_{g_{+}}\pm i\lambda)\sigma^{\pm}=0 \text{ with } x^{\lambda-\frac{n}{2}}\sigma^{\pm}|_{x=0}=\psi^{\pm},$$
if and only if $\Psi^{\pm} =x^{\lambda-\frac{n}{2}}\sigma^{\pm}$ is a solution for
$$\left\{\begin{array}{ll}
D_{\overline{g}}^{f_{1}}\Psi^{\pm}\pm \frac{i\lambda}{x}\Psi^{\pm}=0,\\
\Psi^{\pm}|_{M}=\psi^{\pm}.
\end{array}
\right.
$$
In particular, if $\lambda \in (0,\frac{1}{2})$ then we have $$\D_{h}^{2\lambda}\psi^{\pm}=\frac{d_{\lambda}}{2\lambda}\lim_{x\to 0}x^{1-2\lambda}\partial_{x}\Psi^{\pm}.$$
\end{proposition}

\section{Energy Identities}
We first start by stating the obvious energy identities that follow from the eigenvalue problem on $(X^{n+1},g_{+},\Sigma_{g_{+}}X)$, which is usually referred to as the renormalized energy identity:
\begin{proposition}\label{enbas}
Assume that $\sigma$ satisfies (\ref{eq1}) with $\lambda \in (0,\frac{n}{2})$. If $\lambda \not \in \frac{1}{2}+\N$ then we have
\begin{equation}\label{FP}
F.P.\int_{x>\varepsilon}|\nabla \sigma|^{2}+(\lambda^{2}-\frac{n(n+1)}{4})|\sigma|^{2}\ dv_{g_{+}}=\frac{n}{d_{\lambda}}\int_{M}\langle \D_{h}^{2\lambda}\psi,\psi\rangle\ dv_{h}.
\end{equation}
and if $\lambda\in \frac{1}{2}+\N$, then we have
\begin{equation}
L.P.\int_{x>\varepsilon}|\nabla \sigma|^{2}+(\lambda^{2}-\frac{n(n+1)}{4})|\sigma|^{2}\ dv_{g_{+}}=\frac{n}{\tilde{d}_{\lambda}}\int_{M}\langle \D_{h}^{2\lambda} \psi,\psi \rangle \ dv_{h}.
\end{equation}
In particular, 
\begin{itemize}
\item if $\lambda \in (0,\frac{1}{2})$ then 
$$\int_{x>\varepsilon}|\nabla \sigma|^{2}+(\lambda^{2}-\frac{n(n+1)}{2})|\sigma|^{2}\ dv_{g_{+}}=\frac{(\frac{n}{2}-\lambda)}{\varepsilon^{2\lambda}}\int_{M}|\psi|^{2}\ dv_{h}+\frac{n}{d_{\lambda}}\int_{M} \langle \D_{h}^{2\lambda}\psi,\psi \rangle \ dv_{h} +o(1).$$
\item If $\lambda \in (\frac{1}{2},1)$ we have
\begin{align}
\int_{x>\varepsilon}|\nabla \sigma|^{2}+(\lambda^{2}-\frac{n(n+1)}{2})|\sigma|^{2}\ dv_{g_{+}}&=\frac{(\frac{n}{2}-\lambda)}{\varepsilon^{2\lambda}}\int_{M}|\psi|^{2}\ dv_{h}+\frac{n+1-2\lambda}{(1-2\lambda)\varepsilon^{2\lambda-1}}\int_{M}\langle \nu\cdot D_{h}\psi,\psi \rangle \ dv_{h}\notag\\
&\quad+\frac{n}{d_{\lambda}}\int_{M} \langle \D_{h}^{2\lambda}\psi,\psi \rangle \ dv_{h} +o(1).\notag
\end{align}
\item If $\lambda=\frac{1}{2}$ then 
\begin{align}
\int_{x>\varepsilon}|\nabla \sigma|^{2}+(\lambda^{2}-\frac{n(n+1)}{2})|\sigma|^{2}\ dv_{g_{+}}&=\frac{n-1}{2\varepsilon}\int_{M}|\psi|^{2}\ dv_{h}+\frac{n\ln(\varepsilon)}{\tilde{d}_{\frac{1}{2}}}\int_{M} \langle \nu\cdot D_{h}\psi,\psi \rangle \ dv_{h}\notag\\
&\quad+ \frac{1}{\tilde{d}_{\frac{1}{2}}}\int_{M}\langle \nu\cdot D_{h}\psi,\psi \rangle+o(1).\notag
\end{align}
\end{itemize}
\end{proposition}
{\it Proof:}
This result follows directly from a simple integration by parts. Indeed, we have
\begin{align}
0&=\int_{x>\varepsilon}\langle \Big(-\Delta_{g_{+}}+(\lambda^{2}-\frac{n(n+1)}{4})\Big)\sigma,\sigma\rangle \ dv_{g_{+}}\notag\\
&=\int_{x>\varepsilon} |\nabla \sigma|^{2}+(\lambda^{2}-\frac{n(n+1)}{4})|\sigma|^{2}\ dv_{g_{+}}-\int_{x=\varepsilon}\langle \nabla_{\nu} \sigma, \sigma \rangle \ dv_{h_{\varepsilon}}.\notag
\end{align}
We also recall that $dv_{g_{+}}=x^{-n-1}dv_{h_{x}}dx$ and 
$$dv_{h_{x}}=dv_{h}\left\{\begin{array}{ll}
\sum_{k=0}^{\frac{n}{2}}v^{(2k)}(h)x^{2k} +o(x^{n})\text{ if $n$ is even}.\\
\\
\sum_{k=0}^{N}v^{2k}(h)x^{2k}+o(x^{2N}) \text{ for any $N\in \N$, if $n$ is odd}.
\end{array}
\right.$$
The conclusion then follows from the expansion of $\sigma$ near the boundary in Theorem \ref{thmscat} .
\hfill$\Box$
\begin{remark}
In general for $\lambda <\frac{n}{2}$, one can write
\begin{align}
\int_{x>\varepsilon}|\nabla \sigma|^{2}+\Big(\lambda^{2}-\frac{n(n+1)}{2}\Big)|\sigma|^{2}\ dv_{g_{+}}&=\varepsilon^{-2\lambda}I_{1}(\psi)+\cdots+\varepsilon^{-2\lambda+k}I_{k}(\psi)\notag\\
&\quad+\frac{n}{d_{\lambda}}\int_{M} \langle \D_{h}^{2\lambda}\psi,\psi \rangle \ dv_{h} +o(1),\notag
\end{align}
where $k=[2\lambda]$ amd the $I_{j}$ are integrals with expression involving the renormalized volumes $v^{(2k)}(h)$.
\end{remark}

We aim, next, to have more precise energy estimates that do not blow up at the boundary. For this purpose, instead of working with $\sigma$, we shift to the use of the extension result of Theorem \ref{thmext1}:
\begin{theorem}\label{en1}
Assume that $\Psi$ satisfies $(\ref{BVP})$ then
\begin{itemize}
\item If $\lambda \in (0,\frac{1}{2})$ then
$$\int_{\bar{X}^{n+1}}\Big( |\nabla^{\overline{g}}\Psi|^{2}+\frac{|\nabla^{\overline{g}} x|}{x}\langle [\nabla^{\overline{g}}_{\nu}+\nu\cdot D_{\overline{g}}] \Psi,\Psi\rangle +E_{\lambda}(x)|\Psi|^{2}\Big)x^{1-2\lambda} \ dv_{\bar{g}}=\frac{2\lambda}{d_{\lambda}}\int_{M} \langle \D_{h}^{2\lambda}\psi,\psi\rangle \ dv_{h}.$$
\item If $\lambda \in (\frac{1}{2},1)$ then
\begin{align}
\int_{x>\varepsilon}\Big( |\nabla^{\overline{g}}\Psi|^{2}+\frac{|\nabla^{\tilde{g}}x|}{x}\langle [\nabla^{\overline{g}}_{\nu}+\nu\cdot D_{\overline{g}}] \Psi,\Psi\rangle +E_{\lambda}(x)|\Psi|^{2}\Big)x^{1-2\lambda} \ dv_{\bar{g}}&=\frac{\varepsilon^{1-2\lambda}}{1-2\lambda}\int_{M}\langle\nu \cdot D_{h}\psi,\psi \rangle \ dv_{h}\notag\\
&\quad+\frac{2\lambda}{d_{\lambda}}\int_{M}\langle \D_{h}^{2\lambda}\psi,\psi\rangle \ dv_{h}.\notag
\end{align}

\end{itemize}
\end{theorem}
This is a straightforward conclusion from the boundary value problem (\ref{BVP}) and the integration by parts for the weighted Laplacian as stated in Proposition \ref{int}. In particular, when $M=\R^{n}$ we have
\begin{corollary}
Assume that $\Psi$ satisfies $(\ref{BVP})$ on $\R^{n+1}_{+}$, then for $\lambda \in (0,\frac{1}{2})$ we have
$$\int_{\R^{n+1}_{+}}t^{1-2\lambda}\Big(|\nabla \Psi|^{2}+\frac{1}{t}\langle D_{\R^{n}}\Psi,\Psi \rangle \Big)\ dv_{\R^{n+1}}= \frac{2\lambda}{d_{\lambda}}\int_{\R^{n}}\langle \D_{\R^{n}}^{\lambda}\psi,\psi \rangle \ dv_{\R^{n}},$$
and if $\lambda \in (\frac{1}{2},1)$ we have
\begin{align}
\int_{t>\varepsilon}t^{1-2\lambda}\Big(|\nabla \Psi|^{2}+\frac{1}{t}\langle D_{\R^{n}}\Psi,\Psi \rangle \Big)\ dv_{\R^{n+1}}&= \frac{\varepsilon^{1-2\lambda}}{1-2\lambda}\int_{\R^{n}}\langle \nu\cdot D_{\R^{n}} \psi,\psi\rangle \ dv_{\R^{n}}\notag\\
&\quad+\frac{2\lambda}{d_{\lambda}}\int_{\R^{n}}\langle \D_{\R^{n}}^{\lambda}\psi,\psi \rangle \ dv_{\R^{n}}.\notag
\end{align}
\end{corollary}

\subsection{General Energy Inequalities}
Using a $\lambda$-admissible defining function provides the following result:
\begin{proposition}\label{en2}
Assume that $\lambda \in (0,\frac{1}{2})$ and fix $\psi \in C^{\infty}(M,\Sigma_{h}M)$, then there exists a unique defining function $\rho$ with the expansion $\rho=x\Big(1+\frac{Q_{h}^{\lambda}}{c_{\lambda}}x^{2\lambda}+o(x^{2\lambda})\Big)$ with $c_{\lambda}=2^{2\lambda}\frac{\Gamma(\lambda)}{\Gamma(-\lambda)}$, such that if $\Psi$ satisfies the $(\ref{BVP})$ with the defining function $\rho$, we have

\begin{align}
&\int_{\overline{X}^{n+1}}\Big( |\nabla^{\tilde{g}}\Psi|^{2}+\frac{|\nabla^{\tilde{g}} \rho|}{\rho}\langle \Big[\nabla^{\tilde{g}}_{\nu}+\nu\cdot D_{\tilde{g}}\Big] \Psi,\Psi\rangle +\frac{R_{\tilde{g}}}{4(1-2\lambda)}|\Psi|^{2}\Big)\rho^{1-2\lambda} \ dv_{\tilde{g}}\notag\\
&=\frac{2\lambda}{d_{\lambda}}\int_{M} \langle \D_{h}^{2\lambda}\psi,\psi\rangle \ dv_{h}-\frac{\lambda(n-2\lambda)}{c_{\lambda}}\int_{M} Q_{h}^{\lambda}|\psi|^{2} \ dv_{h},\notag
\end{align}
where $\tilde{g}=\rho^{2}g_{+}$.
\end{proposition}

An interesting identity appears here if we take $\psi=\psi^{\pm}\in C^{\infty}(\Sigma^{\pm}_{h}M)$. Indeed, $\D_{h}^{2\lambda}$ maps $C^{\infty}(\Sigma^{\pm}M)$ to $C^{\infty}(\Sigma^{\mp}M)$. Hence, $\langle \D_{h}^{2\lambda}\psi^{\pm},\psi^{\pm}\rangle =0$ and we get:
$$\int_{\overline{X}^{n+1}}\Big( |\nabla^{\tilde{g}}\Psi|^{2}+\frac{|\nabla^{\tilde{g}} \rho|}{\rho}\langle \Big[\nabla^{\tilde{g}}_{\nu}+\nu\cdot \tilde{D}\Big] \Psi,\Psi\rangle +\frac{\tilde{R}}{4(1-2\lambda)}|\Psi|^{2}\Big)\rho^{1-2\lambda} \ dv_{\tilde{g}}=-\frac{\lambda(n-2\lambda)}{c_{\lambda}}\int_{M} Q_{h}^{\lambda}|\psi^{\pm}|^{2} \ dv_{h}.$$
Using the vanishing of $\langle \D_{h}^{2\lambda}\psi^{\pm},\psi^{\pm}\rangle$, we can improve this previous identity as follow:
\begin{theorem}\label{en3}
Assume that $\lambda \in (0,\frac{1}{4})$ and let $\psi^{\pm}\in C^{\infty}(M,\Sigma^{\pm}M)$, then there exists a $2\lambda$-admissible defining function $\rho$ with the expansion $\rho=x\Big(1+\frac{Q_{h}^{2\lambda}}{c_{2\lambda}}x^{4\lambda}+o(x^{4\lambda})\Big)$ such that if $\Psi$ satisfies $(\ref{BVP})$ with the defining function $\rho$, then
\begin{align}
&\int_{\overline{X}^{n+1}}\Big( |\nabla^{\tilde{g}}\Psi|^{2}+\frac{|\nabla^{\tilde{g}} \rho|}{\rho}\langle \Big[(1-2\lambda)\nabla^{\tilde{g}}_{\nu}+\nu\cdot \tilde{D}\Big] \Psi,\Psi\rangle +E_{\lambda}(\rho)|\Psi|^{2}\Big)\rho^{1-4\lambda} \ dv_{\tilde{g}}\notag\\
&=\frac{2\lambda}{d_{\lambda}^{2}}\int_{M}|\D_{h}^{2\lambda}\psi^{\pm}|^{2}\ dv_{h}-\frac{2\lambda(n-4\lambda)}{c_{2\lambda}}\int_{M} Q_{h}^{2\lambda}|\psi^{\pm}|^{2} \ dv_{h}.\notag
\end{align}

\end{theorem}
{\it Proof:}
The main idea here is to manipulate the PDE of (\ref{BVP}) in order to make it appropriate for the scale $2\lambda$. Indeed, we have
\begin{align}
0&=-\Delta_{\tilde{g}}^{f}\Psi+\frac{|\nabla^{\tilde{g}} \rho|}{\rho}\Big[\nabla^{\tilde{g}}_{\nu}+\nu\cdot D_{\tilde{g}}\Big]\Psi +E_{\lambda}(\rho)\Psi\notag\\
&=-\Delta_{\tilde{g}}^{\tilde{f}}\Psi+\frac{|\nabla^{\tilde{g}} \rho|}{\rho}\Big[(1-2\lambda)\nabla^{\tilde{g}}_{\nu}+\nu\cdot D_{\tilde{g}}\Big]\Psi +E_{\lambda}(\rho)\Psi,\notag
\end{align}
where here $\tilde{f}=\ln(x^{4\lambda-1})$. Now, applying $\Psi$ and integrating by parts with respect to the measure $x^{1-4\lambda}dv_{\tilde{g}}$ leads to:

$$\int_{\rho>\varepsilon}\Big(|\nabla^{\tilde{g}}\Psi|^{2}+\frac{|\nabla^{\tilde{g}} \rho|}{\rho}\langle \Big[(1-2\lambda)\nabla^{\tilde{g}}_{\nu}+\nu\cdot D_{\tilde{g}}\Big] \Psi,\Psi\rangle +E_{\lambda}(\rho)|\Psi|^{2}\Big)\rho^{1-4\lambda} \ dv_{\tilde{g}}=\varepsilon^{1-4\lambda}\int_{\rho=\varepsilon}\langle \nabla^{\tilde{g}}_{\nu}\Psi,\Psi \rangle \ dv_{h_{\varepsilon}}.$$
But
$$\Psi = \psi^{\pm}+S(\lambda)\psi^{\pm}\rho^{2\lambda}-(\frac{n}{2}-\lambda)\frac{Q_{h}^{2\lambda}}{c_{2\lambda}}\psi^{\pm}\rho^{4\lambda}+o(\rho^{4\lambda}).$$
Thus, 
$$\langle \partial_{\rho} \Psi,\Psi\rangle = \Big(2\lambda |S(\lambda)\psi^{\pm}|^{2}-4\lambda(\frac{n}{2}-\lambda)\frac{Q_{h}^{2\lambda}}{c_{2\lambda}}|\psi^{\pm}|^{2}\Big)\rho^{4\lambda-1}+o(\rho^{4\lambda-1}).$$
It follows that,
$$\varepsilon^{1-4\lambda}\lim_{\varepsilon\to 0}\int_{\rho=\varepsilon}\langle \nabla^{\tilde{g}}_{\nu}\Psi,\Psi \rangle \ dv_{h_{\varepsilon}}=\frac{2\lambda}{d_{\lambda}^{2}}\int_{M}|\D_{h}^{2\lambda}\psi^{\pm}|^{2}\ dv_{h}-\frac{4\lambda(n-4\lambda)}{c_{2\lambda}}\int_{M} Q_{h}^{2\lambda}|\psi^{\pm}|^{2} \ dv_{h}.$$
\hfill$\Box$

The last identity shows a link between the operator $(\D_{h}^{2\lambda})^{2}$ and $Q_{h}^{2\lambda}$. This is to be compared to the classical Bochner-Lichnerowicz identity linking $D_{g}^{2}$ and the scalar curvature.\\

\subsection{Sobolev-type Inequality}
We move now to proving a general weighted Sobolev type inequality for spinors for the range $\lambda \in (0,\frac{1}{2})$. Given $\psi \in C^{\infty}(\Sigma_{h}M)$, and fix a $\lambda$-admissible defining function of the form $\rho=x\Big(1+Wx^{2\lambda}+o(x^{2\lambda})\Big)$, we define the space $C^{\lambda}_{\psi}$ by
$$C^{\lambda}_{\psi}:=\Big\{\Psi \in C^{\infty}(\Sigma_{\tilde{g}}X)\cap C^{0}(\Sigma_{\tilde{g}}\overline{X}); \Psi =\psi+\rho^{2\lambda}\varphi +o(\rho^{2\lambda}); \varphi \in C^{\infty}(\Sigma_{h}M)\Big\},$$
and $$C^{\lambda}:=\bigcup_{\psi\in C^{\infty}(\Sigma_{h}M)}C^{\lambda}_{\psi}.$$
We define next the operator
$$K_{\lambda}\Psi:=-\Delta_{\tilde{g}}^{f}\Psi+\frac{|\nabla^{\tilde{g}} \rho|}{\rho}(\nabla^{\tilde{g}}_{\nu}+\nu\cdot D_{\tilde{g}})\Psi +E_{\lambda}(\rho)\Psi,$$
along with the two boundary operators
$$B_{0}^{\lambda}\Psi=\Psi|_{M},$$ and $$B_{2\lambda}^{\lambda}\Psi:=-\lim_{\rho\to 0}\rho^{1-2\lambda}\partial_{\rho}\Psi+\lambda(n-2\lambda)W B_{0}^{\lambda}\Psi.$$
\begin{lemma}
If $\hat{g}=e^{2u}g$ and $\hat{\rho}=e^{u}\rho$ with $u \in C^{\infty}(\overline{X}^{n+1})$, then we have for any $\Psi\in C^{\lambda}$
\begin{itemize}
\item $\widehat{B}_{2\lambda}^{0}(\Psi)=e^{-\frac{n-2\lambda}{2}u|_{M}}B_{2\lambda}^{0}(e^{\frac{n-2\lambda}{2}u}\Psi),$
\item $\widehat{B}_{2\lambda}^{\lambda}(\Psi)=e^{-\frac{n+2\lambda}{2}u|_{M}}B_{2\lambda}^{\lambda}(e^{\frac{n-2\lambda}{2}u}\Psi).$
\end{itemize}
\end{lemma}
{\it Proof:}
The conformal invariance of $B_{2\lambda}^{0}$ is trivial since it is just a restriction to the boundary. So we will focus on the conformal invariance of $B_{2\lambda}^{\lambda}$. By definition,
\begin{align}
B_{2\lambda}^{\lambda}(e^{\frac{n-2\lambda}{2}u}\Psi)&=-\lim_{\rho \to 0}\rho^{1-2\lambda}\partial_{\rho}(e^{\frac{n-2\lambda}{2}u}\Psi)+\lambda(n-2\lambda)WB_{2\lambda}^{0}(e^{\frac{n-2\lambda}{2}u}\Psi)\notag\\
&=e^{\frac{n}{2}u|_{M}}\Big(-\lim_{\hat{\rho}\to 0}e^{\lambda u}\hat{\rho}^{1-2\lambda}\partial_{\hat{\rho}}\Psi + \lambda(n-2\lambda)e^{\lambda u|_{M}}\widehat{W}\widehat{B}_{2\lambda}^{0}(\Psi)\Big)\notag\\
&=e^{\frac{n+2\lambda}{2}u|_{M}}\widehat{B}_{2\lambda}^{\lambda}(\Psi),
\end{align}
where we used the fact that $\widehat{W}=e^{-2\lambda u}W$, which follows from
\begin{align}
\hat{\rho}&=\hat{x}(1+\widehat{W}\hat{x}^{2\lambda}+o(\hat{x}^{2\lambda}))\notag\\
&=e^{u}\rho=e^{u}x(1+We^{-2\lambda u}\hat{x}^{2\lambda}+o(\hat{x}^{2\lambda})).
\end{align}
\hfill $\Box$

In this setting, we define the energy functional $J$ on the weighted Sobolev space $W^{1,2}(\overline{X},\rho^{1-2\lambda}dv_{\tilde{g}})$ by
$$J(\Psi):=\int_{\overline{X}^{n+1}}\langle\Psi, K_{\lambda}\Psi\rangle\  \rho^{1-2\lambda}dv_{\tilde{g}}+\int_{M}\langle B_{0}^{\lambda}\Psi, B_{2\lambda}^{\lambda}\Psi \rangle \ dv_{h}.$$
The corresponding quadratic form is then defined by $$Q(\Psi,\Phi):=\int_{\overline{X}^{n+1}}\langle\Psi, K_{\lambda}\Phi\rangle\  \rho^{1-2\lambda}dv_{\tilde{g}}+\int_{M}\langle B_{0}^{\lambda}\Psi, B_{2\lambda}^{\lambda}\Phi\rangle \ dv_{h}.$$
Hence, we have
\begin{align}
Q(\Psi,\Phi)&=\int_{\overline{X}}\Big(\langle \nabla^{\tilde{g}} \Psi, \nabla^{\tilde{g}} \Phi \rangle +\frac{|\nabla^{\tilde{g}} \rho|}{\rho}\langle (\nabla^{\tilde{g}}_{\nu}+\nu\cdot D_{\tilde{g}})\Phi,\Psi\rangle +E_{\lambda}(\rho)\langle \Psi,\Phi\rangle\Big) \rho^{1-2\lambda}\ dv_{\tilde{g}}\notag\\
&\quad+\lambda(n-2\lambda)\int_{M}W\langle B_{0}^{\lambda} \Psi, B_{0}^{\lambda}\Phi \rangle \ dv_{h},\notag
\end{align}
and 
$$J(\Psi)=\int_{\overline{X}}\Big(|\nabla^{\tilde{g}} \Psi|^{2} +\frac{|\nabla^{\tilde{g}} \rho|}{\rho}\langle (\nabla^{\tilde{g}}_{\nu}+\nu\cdot D_{\tilde{g}})\Psi,\Psi\rangle +E_{\lambda}(\rho)|\Psi|^{2}\Big) \rho^{1-2\lambda}\ dv_{\tilde{g}}+\frac{(n-2\lambda)}{2}\int_{M}W|B_{0}^{\lambda} \Psi|^{2} \ dv_{h}.$$
\begin{proposition}\label{lowerb}
Using the notations above, we have
$$\inf_{\Psi\in C^{\lambda}_{\psi}}J(\Psi)>-\infty$$
\end{proposition}
{\it Proof:}
We consider the first eigenvalue of the operator $K_{\lambda}$ defined by
$$\lambda_{1}(K_{\lambda}):=\inf\{J(\Phi), \Phi \in C^{\lambda}_{0}; \int_{\overline{X}}|\Phi|^{2} \rho^{1-2\lambda}\ dv_{\tilde{g}}=1\}.$$
We claim that $\lambda_{1}(K_{\lambda})>0$. Indeed, if $\sigma=\rho^{\frac{n}{2}-\lambda}\Phi$ with $\Phi \in C^{\lambda}_{0}$, then
\begin{align}
J(\Phi)&=\int_{\overline{X}}\langle (D^{2}_{g_{+}}+\lambda^{2}) \sigma,\sigma \rangle \ dv_{g_{+}}\notag\\
&\geq \lambda^{2}\int_{\overline{X}}|\sigma|^{2}\ dv_{g_{+}}\notag\\
&= \lambda^{2}\int_{\overline{X}}\frac{1}{\rho^{2}}|\Phi|^{2} \rho^{1-2\lambda}\ dv_{\tilde{g}}\notag \\
&\geq \lambda^{2}\int_{\overline{X}}|\Phi|^{2} \rho^{1-2\lambda}\ dv_{\tilde{g}}.\notag
\end{align}
Now, if we fix $\Psi \in C^{\lambda}_{\psi}$ and $\Phi \in C^{\lambda}_{0}$ we have 
\begin{align}
J(\Psi+\Phi)&=J(\Psi)+J(\Phi)+2Q(\Psi,\Phi)\notag\\
&\geq \lambda_{1}(K_{\lambda})\int_{\overline{X}}|\Phi|^{2}\rho^{1-2\lambda}\ dv_{\tilde{g}}+J(\Psi)\notag\\
&\quad -2\Big(\int_{\overline{X}}|K_{\lambda}\Psi|^{2} \rho^{1-2\lambda}\ dv_{\tilde{g}}\Big)^{\frac{1}{2}}\Big(\int_{\overline{X}}|\Phi|^{2} \rho^{1-2\lambda}\ dv_{\tilde{g}}\Big)^{\frac{1}{2}}.
\end{align}
and the conclusion then follows from Young's inequality.
\hfill$\Box$

\begin{theorem}\label{sob}
For all $\Psi \in W^{1,2}(\Sigma_{\tilde{g}}\overline{X},\rho^{1-2\lambda}dv_{\tilde{g}})$ with $Tr(\Psi)=\psi \in C^{\infty}(\Sigma M)$ we have
\begin{align}
\frac{2\lambda}{d_{\lambda}}\int_{M}\langle \D_{h}^{2\lambda}\psi,\psi \rangle\ dv_{h}-\lambda(n-2\lambda)\int_{M}W|\psi|^{2}\ dv_{h}&\leq \int_{\overline{X}}\Big(|\nabla^{\tilde{g}} \Psi|^{2} +\frac{|\nabla^{\tilde{g}} \rho|}{\rho}\langle (\nabla^{\tilde{g}}_{\nu}+\nu\cdot D_{\tilde{g}})\Psi,\Psi\rangle \notag\\
&\qquad+E_{\lambda}(\rho)|\Psi|^{2}\Big) \rho^{1-2\lambda}\ dv_{\tilde{g}}
\end{align}
Moreover, there is equality if and only if $\Psi$ satisfies $K_{\lambda}\Psi=0$.
\end{theorem}
{\it Proof:}
First, notice that $J$ is well defined on the weighted Sobolev space $W^{1,2}(\overline{X},\rho^{1-2\lambda}dv_{\tilde{g}})$. Moreover, as we have seen in Proposition \ref{lowerb}, $J$ is bounded below on $C^{\lambda}_{\psi}$ for any $\psi \in C^{\infty}(M)$. It is clear then that any minimizing sequence $\Psi_{k}\in W^{1,2}(\overline{X},\rho^{1-2\lambda}dv_{\tilde{g}})$ with $Tr(\Psi_{k})=\psi$ has a convergent subsequence to a minimizer $\Psi$ of $J$. This minimizer satisfies the boundary value problem (\ref{BVP}), that is, the equation $K_{\lambda}\Psi=0$ and $\Psi|_{x=0}=\psi$. Hence, we have from Proposition \ref{en2} that
$$J(\Psi)=\frac{2\lambda}{d_{\lambda}} \int_{M}\langle \D_{h}^{2\lambda}\psi,\psi \rangle \ dv_{h},$$
which finishes the proof of the theorem.
\hfill$\Box$

Notice that for the usual pick of a defining function, we can take $W=\frac{Q_{h}^{\lambda}}{c_{\lambda}}$, leading to the inequality
\begin{align}
\frac{2\lambda}{d_{\lambda}} \int_{M}\langle \D_{h}^{2\lambda}\psi,\psi \rangle \ dv_{h}-\frac{\lambda(n-2\lambda)}{c_{\lambda}}\int_{M}Q_{h}^{\lambda}|\psi|^{2}\ dv_{g} &\leq \int_{\overline{X}}\Big(|\nabla^{\tilde{g}} \Psi|^{2} +\frac{|\nabla^{\tilde{g}} \rho|}{\rho}\langle (\nabla^{\tilde{g}}_{\nu}+\nu\cdot D_{\tilde{g}})\Psi,\Psi\rangle\notag\\
&\qquad +\frac{R_{\tilde{g}}}{4(1-2\lambda)}|\Psi|^{2}\Big) \rho^{1-2\lambda}\ dv_{\tilde{g}}.
\end{align}

\section{Critical Operator and $Q$-curvature}
In order to proceed with the critical case $\lambda=\frac{n}{2}$, we need to analyze the solution $\sigma(\lambda)$ of (\ref{eq1}) near $\frac{n}{2}$ when $n$ is odd, where we will be following the ideas in \cite{GZ,T}. For that purpose, let us recall the construction of the solution $\sigma$.
\begin{lemma}\label{lemsig}
Given $\psi \in C^{\infty}(\Sigma_{h}M)$, there exists $\sigma_{\infty}(\lambda)\psi\in x^{\frac{n}{2}-\lambda}C^{\infty}(\Sigma_{\overline{g}}\overline{X})$ such that:
\begin{itemize}
\item[i)] $\frac{1}{\Gamma(\frac{1}{2}-\lambda)}\sigma_{\infty}(\lambda)$ is holomorphic in $\lambda$.
\item[ii)] $\sigma_{\infty}(\lambda)\psi=x^{\frac{n}{2}-\lambda}\psi+o(x^{\frac{n}{2}-\lambda}).$
\end{itemize}
Moreover, 
$$\sigma_{\infty}(\lambda)=x^{\frac{n}{2}-\lambda}\sum_{j=0}^{\infty}x^{k}p_{j,\lambda},$$ 
where the $p_{j,\lambda}$ satisfy the following properties:
\begin{itemize}
\item[a)] $p_{j,\lambda}$ is a differential operator of order $j$ and $p_{0,\lambda}=1$.
\item[b)] $p_{j,\lambda}$ is meromorphic near $\lambda =\frac{n}{2}$ with an isolated simple pole at $\lambda=\frac{n}{2}$ if $j\geq n$.
\item[c)] $(D_{g_{+}}^{2}+\lambda^{2})\Big(x^{\frac{n}{2}-\lambda}\sum_{j=0}^{k}x^{j}p_{j,\lambda}\psi\Big) \in x^{\frac{n}{2}-\lambda+k+1}C^{\infty}(\Sigma_{\overline{g}}\overline{X})$. In particular, $(D_{g_{+}}^{2}+\lambda^{2})\sigma_{\infty}(\lambda)\psi\in x^{\infty}C^{\infty}(\Sigma_{\overline{g}}\overline{X}).$
\end{itemize}
\end{lemma}
The proof of this lemma follows exactly the iterative construction done in \cite{F} and \cite{GMP} for the case of scalar functions.
\begin{lemma}\label{lemres}
Using the notations above, we have
$$Res_{\lambda=\frac{n}{2}}p_{n+j,\lambda}\psi=p_{j,-\lambda}(Res_{\lambda=\frac{n}{2}}p_{n,\lambda})\psi.$$
\end{lemma}
{\it Proof:}
First notice that $\sigma_{\infty}(\lambda)\psi$ and $\sigma_{\infty}(-\lambda)(Res_{\lambda=\frac{n}{2}}p_{n,\lambda})\psi$ satisfy similar equations:
 $$(D_{g_{+}}^{2}+\lambda^{2})\beta=O(x^{\infty}).$$
Taking the residue at $\lambda=\frac{n}{2}$ in Lemma \ref{lemsig}, part $c)$, we see that 
$$(D_{g_{+}}^{2}+\frac{n^{2}}{4})\Big(x^{n}\sum_{j=0}^{k}Res_{\lambda=\frac{n}{2}}p_{n+j,\lambda}\psi\Big)\in x^{n+k+1}C^{\infty}(\Sigma_{\overline{g}}\overline{X}),$$ and similarly
$$(D_{g_{+}}^{2}+\frac{n^{2}}{4})\Big(x^{n}\sum_{j=0}^{k}p_{j,-\lambda}Res_{\lambda=\frac{n}{2}}p_{n,\lambda}\psi\Big)\in x^{n+k+1}C^{\infty}(\Sigma_{\overline{g}}\overline{X}).$$
Therefore, by the iterative construction of $p_{j,-\lambda}$ we get that
$$Res_{\lambda=\frac{n}{2}}p_{n+j,\lambda}\psi=p_{j,-\lambda}(Res_{\lambda=\frac{n}{2}}p_{n,\lambda})\psi.$$
\hfill$\Box$

We are now in position to provide a good expression for the solution $\sigma(\lambda)$:
\begin{theorem}\label{thmres}
Given $\psi \in C^{\infty}(\Sigma_{h}M)$, then there exist families of operators $\sigma_{1}(\lambda), \sigma_{2}(\lambda)$ and $\sigma_{3}(\lambda)\in C^{\infty}(\Sigma_{\overline{g}}\overline{X})$, holomorphic for $\lambda$ in the neighborhood of $\frac{n}{2}$ such that
\begin{itemize}
\item[i)] $\sigma_{1}(\lambda)|_{x=0}=(\lambda-\frac{n}{2})\psi$.
\item[ii)] $\sigma_{1}(\frac{n}{2})+x^{n}\sigma_{2}(\frac{n}{2})\in x^{\infty}C^{\infty}(\Sigma_{\overline{g}}\overline{X})$.
\end{itemize}
Moreover, If $\sigma(\lambda)$ is the unique solution for $(D^{2}_{g_{+}}+\lambda^{2})\sigma(\lambda)=0$ with $x^{-\frac{n}{2}+\lambda}\sigma(\lambda)|_{x=0}=\psi$, then near $\lambda =\frac{n}{2}$ we have:

$$\sigma(\lambda)=\left\{\begin{array}{ll}
  \frac{1}{\lambda-\frac{n}{2}}\Big(x^{\frac{n}{2}-\lambda}\sigma_{1}(\lambda)+x^{\frac{n}{2}+\lambda}\sigma_{2}-\sigma_{1}(\frac{n}{2})-x^{n}\sigma_{2}(\frac{n}{2})\Big)+x^{\frac{n}{2}+\lambda}\sigma_{3}(\lambda) \text{ if $\lambda\not=\frac{n}{2}$ }\\
\\
\dot{\sigma}_{1}(\frac{n}{2})+x^{n}(\dot{\sigma}_{2}(\frac{n}{2})+\sigma_{3}(\frac{n}{2}))+\ln(x)(-\sigma_{1}(\frac{n}{2})+x^{n}\sigma_{2}(\frac{n}{2})) \text{ if $\lambda=\frac{n}{2}$}
\end{array}
\right.
$$
whre $\dot{\sigma}_{i}(\lambda)=\frac{\partial \sigma_{i}}{\partial \lambda}(\lambda)$ for $i=1,2$.
\end{theorem}
{\it Proof:}
We set $\sigma_{1}(\lambda)=x^{\lambda-\frac{n}{2}}(\lambda-\frac{n}{2})\sigma_{\infty}(\lambda)\psi$ and $\sigma_{2}(\lambda)=-x^{-\lambda-\frac{n}{2}} \sigma_{\infty}(-\lambda)Res_{\lambda=\frac{n}{2}}p_{n,\lambda}\psi$. Then clearly, from Lemma \ref{lemres}, $\sigma_{1}$ and $\sigma_{2}$ are holomorphic near $\lambda=\frac{n}{2}$ and satisfy $i)$. If we truncate the expression of $\sigma_{\infty}$ we see that the partial sum of $x^{\frac{n}{2}-\lambda}\sigma_{1}(\lambda)+x^{\frac{n}{2}+\lambda}\sigma_{2}(\lambda)$ can be expressed as
$$x^{\frac{n}{2}-\lambda}\sum_{j=0}^{n+k}(\lambda-\frac{n}{2})x^{j}p_{j,\lambda}\psi -x^{\frac{n}{2}+\lambda}\sum_{j=0}^{k}x^{j}p_{j,-\lambda}Res_{\lambda=\frac{n}{2}}p_{n,\lambda}\psi,$$
which is holomorphic near $\lambda=\frac{n}{2}$ and equals zero at that point and this addresses ii). So we let 
\begin{equation}\label{sigtil}
\tilde{\sigma}_{\infty}(\lambda)=\frac{1}{\lambda -\frac{n}{2}}\Big(x^{\frac{n}{2}-\lambda}\sigma_{1}(\lambda)+x^{\frac{n}{2}+\lambda}\sigma_{2}(\lambda)-\sigma_{1}(\frac{n}{2})-x^{n}\sigma_{2}(\frac{n}{2})\Big),
\end{equation}
and 
$$\sigma(\lambda)=\Big(I-R(\lambda)(D_{g_{+}}^{2}+\lambda^{2})\Big)\tilde{\sigma}_{\infty}(\lambda).$$
Since $(D_{g_{+}}^{2}+\lambda^{2})\tilde{\sigma}_{\infty}(\lambda)=O(x^{\infty})$, we have  $R(\lambda)(D_{g_{+}}^{2}+\lambda^{2})\tilde{\sigma}_{\infty}(\lambda)\in x^{\frac{n}{2}+\lambda}C^{\infty}(\Sigma_{\overline{g}}\overline{X})$. Thus, we define $$\sigma_{3}(\lambda):=-x^{-\lambda-\frac{n}{2}}R(\lambda)(D_{g_{+}}^{2}+\lambda^{2})\tilde{\sigma}_{\infty}(\lambda).$$
Therefore, if $\lambda \not =\frac{n}{2}$ we have
$$(D_{g_{+}}^{2}+\lambda^{2})\sigma(\lambda)=0,$$
and $x^{\lambda-\frac{n}{2}}\sigma(\lambda)|_{x=0}=\psi$. Moreover, $$\sigma(\lambda)=\frac{1}{\lambda -\frac{n}{2}}\Big(x^{\frac{n}{2}-\lambda}\sigma_{1}(\lambda)+x^{\frac{n}{2}+\lambda}\sigma_{2}(\lambda)-\sigma_{1}(\frac{n}{2})-x^{n}\sigma_{2}(\frac{n}{2})\Big)+x^{\frac{n}{2}+\lambda}\sigma_{3}(\lambda),$$
for $\lambda\not=\frac{n}{2}$. Notice now that from $(\ref{sigtil})$, $\sigma$ is holomorphic around $\lambda =\frac{n}{2}$. Therefore, passing to the limit as $\lambda\to \frac{n}{2}$, we get the desired result.
\hfill $\Box$

Based on the result above we see that near $\lambda = \frac{n}{2}$, we have
$$S(\lambda)=\frac{\sigma_{2}(\lambda)|_{x=0}}{\lambda-\frac{n}{2}}+\sigma_{3}(\lambda)|_{x=0}.$$
Thus,
$$\D_{h}^{n}=\tilde{d}_{\frac{n}{2}}\sigma_{2}(\frac{n}{2})|_{x=0}.$$
Hence,
$$\tilde{d}_{\frac{n}{2}}S(\lambda)\psi=\frac{\D_{h}^{n}\psi}{\lambda-\frac{n}{2}}+\tilde{d}_{\frac{n}{2}}\Big(\dot{\sigma}_{2}(\frac{n}{2})+\sigma_{3}(\frac{n}{2})\Big)+O(\lambda-\frac{n}{2}).$$
\begin{definition}\label{defQ}
Using the notations of Theorem \ref{thmres}, if $n$ is odd we define the spinorial $Q$-curvature operator $\Q_{h}:C^{\infty}(\Sigma_{h}M)\to C^{\infty}(\Sigma_{h}M)$ by
$$\Q_{h}(\psi)=-\tilde{d}_{\frac{n}{2}}\Big(\dot{\sigma}_{2}(\frac{n}{2})+\sigma_{3}(\frac{n}{2})\Big).$$
In particular, 
$$\tilde{d}_{\frac{n}{2}}S(\lambda)\psi=\frac{\D_{h}^{n}\psi}{\lambda-\frac{n}{2}}-\Q_{h}(\psi)+O(\lambda-\frac{n}{2}).$$
\end{definition}

We recall that the operator $\D_{h}^{n}$ satisfies the following conformal invariance property:
$$\D^{n}_{\hat{h}}=e^{-nw}\D_{h}^{n}.$$
Hence, $\mathcal{H}:=\ker \D_{h}^{\frac{n}{2}}$ is conformally invariant. We can then deduce the following: 

\begin{theorem}\label{thmQ}
The operator $\Q_{h}: \mathcal{H} \to C^{\infty}(\Sigma_{\bar{g}}\overline{X}^{n+1}|_{M})$ is essentially self-adjoint with the property that for all $\psi \in \mathcal{H}$,
\begin{equation}\label{Qvar}
e^{nw}\Q_{\hat{h}}(\psi)=\Q_{h}(\psi)+\D^{\frac{n}{2}}_{h}(w\psi).
\end{equation}
Moreover, for $\psi_{1}, \psi_{2} \in \mathcal{H}$, we let $Q_{h}^{\psi_{1},\psi_{2}}:=\langle \Q_{h}(\psi_{1}),\psi_{2}\rangle$ and define the operator $\mathcal{L}_{h}^{\psi_{1},\psi_{2}}(\cdot):=\langle \D_{h}^{\frac{n}{2}}(\cdot \psi_{1}),\psi_{2}\rangle$, then we have
$$Q_{\hat{h}}^{\psi_{1},\psi_{2}}e^{nw}=Q_{h}^{\psi_{1},\psi_{2}}+\mathcal{L}_{h}(\omega).$$
In addition, $\overline{Q}^{\psi_{1},\psi_{2}}:=\int_{M}Q_{h}^{\psi_{1},\psi_{2}}\ dv_{h}$ is a conformal invariant.
\end{theorem}
{\it Proof:}
We first start by fixing $\psi \in \mathcal{H}$. Then we set
$$ \Q_{h}(\psi)=-\frac{\partial \D_{h}^{2\lambda}\psi}{\partial \lambda}_{|\lambda =\frac{n}{2}}.$$ 
This definition extends also to the case $n$ is even. Notice that if $n$ is odd, this definition coincides with the one in Definition \ref{defQ}. Indeed,
$$\D^{2\lambda}_{h}=d_{\lambda}S(\lambda)=\frac{\tilde{d}_{\lambda}}{\tilde{d}_{\frac{n}{2}}}\Big(\D_{h}^{n}\psi-(\lambda-\frac{n}{2})\Q_{h}(\psi)+O((\lambda-\frac{n}{2})^{2})\Big).$$
Hence, since $\D_{h}^{n}\psi=0$, we see that
 $$-\frac{\partial \D_{h}^{2\lambda}\psi}{\partial \lambda}_{|\lambda =\frac{n}{2}}=\Q_{h}(\psi).$$
When $n$ is even, $\Q_{h}(\psi)$ is clearly well defined since $\D^{2\lambda}_{h}$ is holomorphic near $\lambda=\frac{n}{2}$.
Now, given $\lambda$ close to $\frac{n}{2}$ and $\hat{h}=e^{2w}h$ with $w\in C^{\infty}(M)$, we have
$$\D_{\hat{h}}^{2\lambda}\psi = e^{-(\frac{n}{2}+\lambda)w}\D_{h}^{2\lambda}e^{(\frac{n}{2}-\lambda)w}\psi.$$
Differentiating this last equality with respect to $\lambda$ and setting $\lambda=\frac{n}{2}$, we get
$$e^{nw}\Q_{\hat{h}}(\psi)=w\D_{h}^{n}\psi+\Q_{h}(\psi)+\D_{h}^{n}(w\psi).$$
Therefore, since $\psi \in \mathcal{H}$, equation (\ref{Qvar}) is satisfied.\\
Next, we let $\psi_{1}, \psi_{2}\in \mathcal{H}$, then clearly, using the identity (\ref{Qvar}), and multiplying by $\psi_{2}$ from the right, we get
$$Q_{\hat{h}}^{\psi_{1},\psi_{2}}e^{nw}=Q_{h}^{\psi_{1},\psi_{2}}+\mathcal{L}_{h}(w).$$
Notice that $\mathcal{L}_{\hat{h}}=e^{-nw}\mathcal{L}_{h}$ and
$$\int_{M}Q_{\hat{h}}^{\psi_{1},\psi_{2}}\ dv_{\hat{h}}=\int_{M}e^{nw}Q_{\hat{h}}^{\psi_{1},\psi_{2}} \ dv_{h}=\int_{M} Q_{h}^{\psi_{1},\psi_{2}}+\mathcal{L}_{h}w \ dv_{h}.$$
But since $\D_{h}^{n}$ is essentially self-adjoint and $\psi_{2}\in \mathcal{H}$, we have
$$\int_{M}\mathcal{L}_{h}(w)\ dv_{h}=\int_{M} \langle \D_{h}^{n}(w\psi_{1}),\psi_{2}\rangle \ dv_{h}=\int_{M} \langle w\psi_{1},\D_{h}^{n}\psi_{2}\rangle \ dv_{h}=0.$$
Thus, $$\int_{M}Q_{\hat{h}}^{\psi_{1},\psi_{2}}\ dv_{\hat{h}}= \int_{M}Q_{h}^{\psi_{1},\psi_{2}}\ dv_{h}.$$
\hfill$\Box$

Using the energy identity (\ref{FP}), we see that if $\psi \in \mathcal{H}$ and $\sigma$ is a solution of (\ref{eq1}), with $x^{\lambda-\frac{n}{2}}\sigma|_{x=0}=\psi$ then for $\lambda$ near $\frac{n}{2}$ and $\lambda \not=\frac{n}{2}$ we have
$$F.P.\int_{x>\varepsilon}|\nabla \sigma|^{2}+(\lambda^{2}-\frac{n(n+1)}{4})|\sigma|^{2}\ dv_{g_{+}}=
-\frac{n}{\tilde{d}_{\frac{n}{2}}}\int_{M}\langle \Q_{h}\psi,\psi\rangle\ dv_{h}+O(\lambda -\frac{n}{2}).$$
\begin{remark}
The spinorial Q-curvature introduced in the previous theorem is similar in nature to the Q-curvature operators acting on forms and studied in \cite{BG,BG2,G,AG}. The $Q$-curvature operators on forms appear naturally when dealing with harmonic forms on the extension. These operators were defined on even dimensional manifolds and they act on closed forms. That is because in that case, closed forms are in the kernel of the corresponding operator $\D_{h}^{n}$ on forms. We also point out that another type of $Q$-curvature was introduced in \cite{GPS} for the Yang-Mills functional by applying the $Q$-curvature operator to the curvature of a Yang-Mills connection (which is a closed form by construction).
\end{remark}

We provide now a spinorial version of the special defining function introduced in \cite{FG1}:
\begin{theorem}\label{thmQ2}
We fix $\psi \in \mathcal{H}$ and let $\sigma(\lambda)$ be the solution of (\ref{eq1}) for $\frac{n}{2}-\varepsilon<\lambda\leq\frac{n}{2}$, with $x^{\lambda-\frac{n}{2}}\sigma(\lambda)|_{x=0}=\psi$, as expressed in Theorem \ref{thmres}. Then, there exists $U$ solving
$$(D_{g_{+}}^{2}+\frac{n^{2}}{4})U=n\sigma(\frac{n}{2})+O(x^{\infty}).$$
Such that:
\begin{itemize}
\item[i)] If $n$ is odd, then there exist $A$, $B$ and $C$ in $C^{\infty}(\Sigma_{\overline{g}}\overline{X})$ such that
$$U=\sigma(\frac{n}{2})\ln(x)+A+x^{n}\ln(x)B+x^{n}\ln^{2}(x)C.$$
Moreover, $A$ is unique modulo $O(x^{n})$, $B$ and $C$ are unique modulo $O(x^{\infty})$ and
$$A|_{x=0}=0, \quad B|_{x=0}=\frac{2}{\tilde{d}_{\frac{n}{2}}}\Q_{h}(\psi), \quad \text{ and } C|_{x=0}=\frac{-2}{\tilde{d}_{\frac{n}{2}}}\D_{h}^{n}\psi.$$
\item[ii)] If $n$ is even, then there exist $\tilde{A}$, $\tilde{B}$ and $\tilde{C}$ in $C^{\infty}(\Sigma_{\overline{g}}\overline{X})$ such that 
$$U=\sigma(\frac{n}{2})\ln(x)+\tilde{A}+x^{n}\ln(x)\tilde{B}+x^{n}\tilde{C},$$
with
$$\tilde{A}|_{x=0}=0,\quad \tilde{B}|_{x=0}=-\frac{1}{d_{\frac{n}{2}}}\D_{h}^{n}\psi, \quad \text{ and } \quad \tilde{C}|_{x=0}=\frac{1}{d_{\frac{n}{2}}}\Q_{h}(\psi).$$
\end{itemize}
\end{theorem}
{\it Proof:}
The construction and uniqueness of such a $U$ follows exactly the same iterative process as in \cite{GZ,F}, starting from the first $U_{0}=n\ln(x)\sigma(\frac{n}{2})$. We start then with the case when $n$ is odd. Recall that $\sigma(\lambda)$ satisfies 
$$(D^{2}_{g_{+}}+\lambda^{2})\sigma(\lambda)=0.$$
Differentiating this last equality and setting $\lambda=\frac{n}{2}$ we get
$$(D^{2}_{g_{+}}+\frac{n^{2}}{4})\dot{\sigma}(\frac{n}{2})=-n\sigma(\frac{n}{2}).$$
We set then $U=-\dot{\sigma}(\frac{n}{2})$. On the other hand, 
\begin{align}
\dot{\sigma}(\frac{n}{2})&=\frac{1}{2}\frac{\partial^{2}}{\partial \lambda^{2}}_{|\lambda=\frac{n}{2}}[x^{\frac{n}{2}-\lambda}\sigma_{1}(\lambda)+x^{\frac{n}{2}+\lambda}\sigma_{2}(\lambda)]+\frac{\partial  x^{\frac{n}{2}+\lambda}\sigma_{3}(\lambda)}{\partial \lambda}_{|\lambda=\frac{n}{2}}\\
&=-\ln(x)\sigma(\frac{n}{2})+2x^{n}\ln(x)[\dot{\sigma}_{2}(\frac{n}{2})+\sigma_{3}(\frac{n}{2})]+2x^{n}\ln^{2}(x)\sigma_{2}(\frac{n}{2})\notag \\
&\quad -\frac{1}{2}\ln^{2}(x)[\sigma_{1}(\frac{n}{2})+x^{n}\sigma_{2}(\frac{n}{2})]+\frac{1}{2}[\ddot{\sigma}_{1}(\frac{n}{2})+x^{n}\ddot{\sigma}_{2}(\frac{n}{2})]+x^{n}\dot{\sigma}_{3}(\frac{n}{2}).\notag
\end{align}
Thus, we can take $$-A=x^{n}\dot{\sigma}_{3}(\frac{n}{2})-\frac{1}{2}\ln^{2}(x)[\sigma_{1}(\frac{n}{2})+x^{n}\sigma_{2}(\frac{n}{2})]+\frac{1}{2}[\ddot{\sigma}_{1}(\frac{n}{2})+x^{n}\ddot{\sigma}_{2}(\frac{n}{2})],$$
$$-B=2[\dot{\sigma}_{2}(\frac{n}{2})+\sigma_{3}(\frac{n}{2})]\quad \text{ and }\quad-C=2\sigma_{2}(\frac{n}{2}).$$
Notice that since $\sigma_{1}(\frac{n}{2})+x^{n}\sigma_{2}(\frac{n}{2})\in x^{\infty}C^{\infty}(\Sigma_{\overline{g}}\overline{X})$ we have that $A\in C^{\infty}(\Sigma_{\overline{g}}\overline{X})$ and $A|_{x=0}=0$. Also,
$$B|_{x=0}=-2[\dot{\sigma}_{2}(\frac{n}{2})+\sigma_{3}(\frac{n}{2})]|_{x=0}=\frac{2}{\tilde{d}_{\frac{n}{2}}}\Q_{h}(\psi),$$
and
$$C|_{x=0}=-2\sigma_{2}(\frac{n}{2})|_{x=0}=-\frac{2}{\tilde{d}_{\frac{n}{2}}}\D_{h}^{n}\psi.$$
In the case where $n$ is even, we have $\sigma(\lambda)=x^{\frac{n}{2}-\lambda}\sigma_{+}(\lambda)+x^{\frac{n}{2}+\lambda}\sigma_{-}(\lambda)$. Moreover, $\sigma_{+}$ and $\sigma_{-}$ are holomorphic near $\lambda=\frac{n}{2}$. Hence, $U=-\dot{\sigma}(\frac{n}{2})$ satisfies
$$(D^{2}_{g_{+}}+\frac{n^{2}}{4})U=n\sigma(\frac{n}{2}).$$
Differentiating the expression of $\sigma$ with respect to $\lambda$ yields
$$U=\ln(x)\sigma_{+}(\frac{n}{2})-\dot{\sigma}_{+}(\frac{n}{2})-x^{n}\ln(x)\sigma_{-}(\frac{n}{2})-x^{n}\dot{\sigma}_{-}(\frac{n}{2}).$$
We let then 
$$\tilde{A}=-\dot{\sigma}_{+}(\frac{n}{2}),\quad \tilde{B}=-\sigma_{-}(\frac{n}{2})\quad \text{ and }\quad \tilde{C}=-\dot{\sigma}_{-}(\frac{n}{2}).$$
In particular, $\tilde{A}|_{x=0}=0$, $\tilde{B}|_{x=0}=-\frac{1}{d_{\frac{n}{2}}}\D_{h}^{n}\psi$ and $\tilde{C}|_{x=0}=\frac{1}{d_{\frac{n}{2}}}\Q_{h}(\psi)$.
\hfill$\Box$

In general, when $\psi \not \in \mathcal{H}$, $\Q_{h}$ follows the following change, under conformal deformation:
$$e^{nw}\Q_{\hat{h}}(\psi)=\Q_{h}(\psi)+\D_{h}^{n}(w\psi)+w\D_{h}^{n}(\psi).$$
This reminds us of the CR case, where $\D_{h}^{n}$ plays the role of the Paneitz operator$P_{\theta}$ and $\Q_{h}$ plays the role of $P'_{\theta}$. This leads us to refine the expansion of $S(\lambda)$ by defining another curvature $\Q'_{h}$ by setting
\begin{equation}\label{Qprime}
\tilde{d}_{\frac{n}{2}}S(\lambda)\psi=\frac{\D_{h}^{n}\psi}{\lambda-\frac{n}{2}}-\Q_{h}(\psi)+(\lambda-\frac{n}{2})\Q'_{h}(\psi)+O((\lambda-\frac{n}{2})^{2}).
\end{equation}
With this definition, we see by looking at the $(\lambda -\frac{n}{2})$ term of (\ref{Qprime}), that under a conformal change $\hat{h}=e^{2w}h$ we have
$$e^{nw}\Big(\frac{1}{2}w^{2}\D_{\hat{h}}(\psi)-w\Q_{\hat{h}}(\psi)+\Q'_{h}(\psi)\Big)=-\frac{1}{2}\D_{h}^{n}(w^{2}\psi)+\Q_{h}(w\psi)+\Q'_{h}(\psi).$$
As is the case in the CR-setting, the $\Q'_{h}$-curvature is relevant when the $\Q_{h}$-curvature vanishes. Hence, we have the following result:
\begin{corollary}
Given $\psi_{0}\in \mathcal{H}\cap \ker \Q_{h}$ then under the conformal change $\hat{h}=e^{2w}h$ we have
$$\Q_{h}(w\psi_{0})+\Q'_{h}(\psi_{0})=e^{nw}\Q'_{\hat{h}}(\psi_{0})-w\D_{h}^{n}(w\psi_{0})+\frac{1}{2}\D_{h}^{n}(w^{2}\psi_{0}).$$
In particular, if $\mathcal{K}:=\ker \D_{h}^{n}(\cdot \psi_{0})\subset C^{\infty}(M)$, we have that for $w\in \mathcal{K}$
$$\Q_{h}(w\psi_{0})+\Q'_{h}(\psi_{0})=e^{nw}\Q'_{\hat{h}}(\psi_{0}) \text{ mod } \mathcal{K}^{\perp}.$$
\end{corollary}
The $\Q'$-curvature is useful to understand the geometry of $M$ only when $\Q_{h}(\psi)=0$. Notice that in contrast with the CR case, the kernel of $\D_{h}^{n}$ is finite dimensional, which imposes strict limitations on the space $\mathcal{K}$ (which can be just $\R$). The space $\mathcal{K}$ corresponds to the set of pluriharmonic functions in the CR setting, which is infinite dimensional.

\appendix

\section{Appendix}
\subsection{The case of the flat space $\R^{n}$}
We consider the Poincar\'{e} half-space model for the hyperbolic space. Namely $\R^{n+1}_{+}$ with the metric $g_{+}=\frac{dt^{2}+dx^{2}}{t^{2}}$. Then the equation 
$(D_{g_{+}}^{2}+\lambda^{2})\sigma=0$ is equivalent to:

$$(\lambda^{2}-(\frac{n}{2})^{2})\sigma+(n-1)t\partial_{t}\sigma+t\partial_{t}\cdot D_{\R^{n}}\sigma+t^{2}D_{\R^{n}}^{2}\sigma-t^{2}\partial_{t}^{2}\sigma=0.$$
We propose to find the expression of the scattering operator $S(\lambda)$ and we will focus on the case $\lambda \in (0,1)$. Hence, after using a Fourier transform in $\R^{n}$ and separating the variables in $\mathcal{F}(\sigma)=f(t)\mathcal{F}(\psi)$ we have that $f$ satisfies:
$$t^{2}f''-(n-1)tf'+[At^{2}+Bt+C]f=0,$$
where $A=\mathcal{F}(\Delta)=-|\xi|^{2}$, $B=-\mathcal{F}(\nu\cdot D_{\R^{n}})=-i\nu\cdot \xi \cdot$ and $C=(\frac{n}{2})^{2}-\lambda^{2}$. Notice that $A$, $B$ and $C$ commute. In fact, one could also use the spectral measure of the operator $\nu\cdot D_{\R^{n}}$ instead of the Fourier transform to reach this equation.\\
The indicial roots correspond to $\gamma_{1}=\frac{n}{2}-\lambda$ and $\gamma_{2}=\frac{n}{2}+\lambda$. We look for a solution of the form $f=t^{\frac{n}{2}+\lambda}y$ with $\lim_{t\to 0}t^{\lambda-\frac{n}{2}}f (t)=1$. This leads to $y$ satisfying
$$ty''+(1+2\lambda)t'+(At+B)y=0.$$
Thus, if we set $y_{1}=e^{-|\xi|t}w(2|\xi|t)$, we see that $w$ satisfies the equation
$$tw''+(1+2\lambda-t)w'-(\frac{1+2\lambda}{2}+\frac{-B}{2|\xi|})w=0,$$
which is Kummer's equation. We recall here some of the properties of this equation:
\begin{lemma}[\cite{RGM}]
Consider the equation $tw''+(b-t)w'-aw=0$, with $a\geq 0$ and $b\geq 1$. Then it has two linearly independent solutions $M(a,b,\cdot)$ and $V(a,b,\cdot)$. When $t\to 0$ we have
$$M(a,b,t)=1+O(t)\quad \text{ and } \quad t^{b-1}V(a,b,t)=\frac{\Gamma(b-1)}{\Gamma(a)}+o(1).$$
Moreover,
$$V(a,b,t)=\frac{\pi}{\sin(\pi b)}\Big(\frac{M(a,b,t)}{\Gamma(1+a-b)\Gamma(b)}-t^{1-b}\frac{M(1+a-b,2-b,t)}{\Gamma(a)\Gamma(2-b)}\Big).$$
\end{lemma}
Therefore, we have
$$f(t)=t^{\frac{n}{2}+\lambda}e^{-|\xi|t}\frac{\Gamma(\frac{-B}{2|\xi|}+\frac{1}{2}+\lambda)}{\Gamma(2\lambda)}(2|\xi|)^{2\lambda}V\Big(\frac{-B}{2|\xi|}+\frac{1}{2}+\lambda,2\lambda+1,2|\xi|t\Big).$$

Thus, near $t=0$ we have
$$f(t)=t^{\frac{n}{2}-\lambda}(1+o(1))+t^{\frac{n}{2}+\lambda}\frac{\Gamma(\frac{-B}{2|\xi|}+\frac{1}{2}+\lambda)\Gamma(-2\lambda)}{\Gamma(\frac{-B}{2|\xi|}+\frac{1}{2}-\lambda)\Gamma(2\lambda)}(2|\xi|)^{2\lambda}(1+o(1)).$$
Therefore, the scattering operator takes the form
$$S(\lambda)=-\frac{\Gamma(-2\lambda)}{\Gamma(2\lambda)}\frac{\Gamma(\frac{\nu\cdot D_{\R^{n}}}{2|D_{\R^{n}}|}+\frac{1}{2}+\lambda)}{\Gamma(\frac{\nu\cdot D_{\R^{n}}}{2|D_{\R^{n}}|}+\frac{1}{2}-\lambda)}2^{2\lambda}|D_{\R^{n}}|^{2\lambda}.$$
We use now the doubling property of the $\Gamma$-function, namely
\begin{equation}\label{doub}
\Gamma(2z)=(2\pi)^{-\frac{1}{2}}2^{2z-\frac{1}{2}}\Gamma(z)\Gamma(z+\frac{1}{2}).
\end{equation}
Hence,
$$S(\lambda)=-2^{-2\lambda}\frac{\Gamma(\frac{1}{2}-\lambda)}{\Gamma(\frac{1}{2}+\lambda)}\frac{\Gamma(-\lambda)}{\Gamma(\lambda)}\frac{\Gamma(\frac{\nu\cdot D_{\R^{n}}}{2|D_{\R^{n}}|}+\frac{1}{2}+\lambda)}{\Gamma(\frac{\nu\cdot D_{\R^{n}}}{2|D_{\R^{n}}|}+\frac{1}{2}-\lambda)}|D_{\R^{n}}|^{2\lambda}.$$
Thus,
$$\D_{\R^{n}}^{2\lambda}=d_{\lambda}S(\lambda)=-\frac{\Gamma(-\lambda)}{\Gamma(\lambda)}\frac{\Gamma(\frac{\nu\cdot D_{\R^{n}}}{2|D_{\R^{n}}|}+\frac{1}{2}+\lambda)}{\Gamma(\frac{\nu\cdot D_{\R^{n}}}{2|D_{\R^{n}}|}+\frac{1}{2}-\lambda)}|D_{\R^{n}}|^{2\lambda}.$$
To conclude, we let $dE_{\ell}$ the spectral measure of the operator $\nu\cdot D_{\R^{n}}$. Then we have,
\begin{align}
\D_{\R^{n}}^{2\lambda}&=-\frac{\Gamma(-\lambda)}{\Gamma(\lambda)}\int_{-\infty}^{+\infty}\frac{\Gamma(\frac{\ell}{2|\ell|}+\frac{1}{2}+\lambda)}{\Gamma(\frac{\ell}{2|\ell|}+\frac{1}{2}-\lambda)}|\ell|^{2\lambda}\ dE_{\ell}\notag\\
&=-\frac{\Gamma(-\lambda)}{\Gamma(\lambda)}\int_{0}^{+\infty}\frac{\Gamma(1+\lambda)}{\Gamma(1-\lambda)}|\ell|^{2\lambda}\ dE_{\ell}-\frac{\Gamma(-\lambda)}{\Gamma(\lambda)}\int_{-\infty}^{0}\frac{\Gamma(\lambda)}{\Gamma(-\lambda)}|\ell|^{2\lambda}\ dE_{\ell}\notag\\
&=\int_{-\infty}^{+\infty}|\ell|^{2\lambda-1}\ell \ dE_{\ell}.\notag
\end{align}
Hence,
$$\D_{\R^{n}}^{2\lambda}=|D_{\R^{n}}|^{2\lambda-1}\nu\cdot D_{\R^{n}}.$$

\subsection{The case of the round sphere}
In order to find the expression of $\D^{2\lambda}_{h_{0}}$ on the standard sphere $(S^{n},h_{0})$ we start from the AdS model of the hyperbolic metric. That is, we consider the manifold $\R\times S^{n}$ with the metric $g_{+}=dy^{2}+\sinh^{2}(y)h_{0}$. For simplicity, we will focus on the case when $n$ is even, but the odd case can be obtained similarly. Fortunately, Camporesi-Higuchi \cite{CH} already computed the expression of the eigenspinors of $D_{g_{+}}$. We will use their approach and their findings to determine the expression of $\D^{2\lambda}_{g_{0}}$.\\
We first recall that in this metric $g_{+}$, the equation $(D^{2}_{g_{+}}+\lambda^{2})\sigma=0$ can be expressed as:
$$\Big[(\partial_{y}+\frac{n}{2}\coth(y))^{2}-\frac{1}{\sinh^{2}(y)}D_{h_{0}}^{2}+\frac{\cosh(y)}{\sinh^{2}(y)}\nu\cdot D_{h_{0}}-\lambda^{2}\Big]\sigma=0.$$
As it was pointed out in \cite{CH}, $\nu\cdot D_{h_{0}}$ commutes with $D_{h_{0}}^{2}$, since $\Big(\nu\cdot D_{h_{0}} \Big)^{2}=D_{h_{0}}^{2}$. Thus, they can be both diagonalized in a common basis of eigenspinors of $\nu\cdot D_{h_{0}}$. We let $(\varphi_{k})_{k\in \Z}$ be the basis of eigenspinors of $D_{h_{0}}$, with the convention that $\varphi_{k}$ corresponds to to the positive eigenvalues when $k>0$, that is for the eigenvalue $\mu_{k}=\frac{n}{2}+k-1$. We will also set $\varphi_{-k}=\nu\cdot \varphi_{k}$ with eigenvalue $-\mu_{k}$. Then we can have an explicit formula for the eigenspinors of $\nu\cdot D_{h_{0}}$. Indeed, for $k>0$, let $\hat{\varphi}_{k,+}=(1-\nu\cdot)\varphi_{-k}$ and $\hat{\varphi}_{k,-}=\nu \hat{\varphi}_{k,+}$, then we have
$$\nu\cdot D_{h_{0}}\hat{\varphi}_{k,+}=\mu_{k}\hat{\varphi}_{k,+},$$
and
$$\nu \cdot D_{h_{0}}\hat{\varphi}_{k,-}=-\mu_{k} \hat{\varphi}_{k,-}.$$
Therefore, it is natural to looks for a solution of the form $\sigma=\sum_{k=1}^{\infty} f_{k}(y)\hat{\varphi}_{k,+}+g_{k}(y)\hat{\varphi}_{k,-}$. These computations were done in \cite{CH} and there, the authors showed that up to a multiplicative constant, one has
$$f_{k}(y)=\sinh^{k-1}(\frac{y}{2})\cosh^{k}(\frac{y}{2}) _{2}F_{1}(\mu_{k}+\frac{1}{2}-\lambda,\mu_{k}+\frac{1}{2}+\lambda,\mu_{k}+\frac{1}{2},-\sinh^{2}(\frac{y}{2})),$$
and
$$g_{k}(y)=\sinh^{k}(\frac{y}{2})\cosh^{k-1}(\frac{y}{2}) _{2}F_{1}(\mu_{k}+\frac{1}{2}-\lambda,\mu_{k}+\frac{1}{2}+\lambda,\mu_{k}+\frac{3}{2},-\sinh^{2}(\frac{y}{2})).$$
Notice that these two solutions are regular at $y=0$, but one needs to analyze their behaviour at infinity, since that is what corresponds to the conformal boundary at infinity. For that purpose, we use the change of variable $y=-\ln(\frac{r}{2})$. Indeed, for $f_{k}$ we have
$$f_{k}(r)=(1-\frac{r^{2}}{4})^{k-1}(1+\frac{r}{2})\sqrt{2r}^{-(2k-1)} {}_{2}F_{1}(\mu_{k}+\frac{1}{2}-\lambda,\mu_{k}+\frac{1}{2}+\lambda,\mu_{k}+\frac{1}{2},-\frac{(1-\frac{r}{2})^{2}}{2r}).$$
We use now the following inversion property for hypergeometric functions: For $z\not \in (0,1)$ and $a-b\not \in \Z$,
\begin{align}
 _{2}F_{1}(a,b,c,z)=&\frac{\Gamma(b-a)\Gamma(c)}{\Gamma(b)\Gamma(c-a)}(-z)^{-a} {}_{2}F_{1}(a,a-c+1,a-b+1,\frac{1}{z})\notag\\
&+\frac{\Gamma(a-b)\Gamma(c)}{\Gamma(a)\Gamma(c-b)}(-z)^{-b} {}_{2}F_{1}(b,b-c+1,b-a+1.\frac{1}{z}).\notag
\end{align}
Thus, 
$$f_{k}(r)=(2r)^{\frac{n}{2}-\lambda}\frac{\Gamma(2\lambda)\Gamma(\mu_{k}+\frac{1}{2})}{\Gamma(\mu_{k}+\frac{1}{2}+\lambda)\Gamma(\lambda)}H_{1}(r)+(2r)^{\frac{n}{2}+\lambda}\frac{\Gamma(-2\lambda)\Gamma(\mu_{k}+\frac{1}{2})}{\Gamma(\mu_{k}+\frac{1}{2}-\lambda)\Gamma(-\lambda)}H_{2}(r),$$
where $$H_{1}(r)=(1-\frac{r^{2}}{4})^{k}(1-\frac{r}{2})^{-2\mu_{k}-2+\lambda}{}_{2}F_{1}\Big(\mu_{k}+\frac{1}{2}-\lambda,1-\lambda,1-2\lambda, \frac{2r}{(1-\frac{r}{2})^{2}}\Big),$$
and $$H_{2}(r)=(1-\frac{r^{2}}{4})^{k}(1-\frac{r}{2})^{-2\mu_{k}-2-2\lambda}{}_{2}F_{1}\Big(\mu_{k}+\frac{1}{2}+\lambda,1+\lambda,1+2\lambda, \frac{2r}{(1-\frac{r}{2})^{2}}\Big).$$
Since $H_{1}(0)=H_{2}(0)=1$ and we are looking for a spectral multiplier of the form $cf_{k}(r)\hat{\varphi}_{k,+}$, we see that the right choice for $c$ that provides the scattering operator operator $S(\lambda)$ is 
$$c=\Big(2^{\frac{n}{2}-\lambda}\frac{\Gamma(2\lambda)\Gamma(\mu_{k}+\frac{1}{2})}{\Gamma(\mu_{k}+\frac{1}{2}+\lambda)\Gamma(\lambda)}\Big)^{-1}.$$
Therefore,
$$S(\lambda)\hat{\varphi}_{k,+}=2^{2\lambda}\frac{\Gamma(\lambda)\Gamma(\mu_{k}+\frac{1}{2}+\lambda)\Gamma(-2\lambda)}{\Gamma(\mu_{k}+\frac{1}{2}-\lambda)\Gamma(2\lambda)\Gamma(-\lambda)}\hat{\varphi}_{k,+}.$$
Using again the doubling identity of the $\Gamma$-function (\ref{doub}), yields:
$$S(\lambda)\hat{\varphi}_{k,+}=2^{-2\lambda}\frac{\Gamma(\mu_{k}+\frac{1}{2}+\lambda)\Gamma(\frac{1}{2}-\lambda)}{\Gamma(\mu_{k}+\frac{1}{2}-\lambda)\Gamma(\frac{1}{2}+\lambda)}\hat{\varphi}_{k,+}.$$
In particular, $$\D_{h_{0}}^{2\lambda}\hat{\varphi}_{k,+}=d_{\lambda}S(\lambda)\hat{\varphi}_{k,+}=\frac{\Gamma(\mu_{k}+\frac{1}{2}+\lambda)}{\Gamma(\mu_{k}+\frac{1}{2}-\lambda)}\hat{\varphi}_{k,+}.$$
A similar process for $g_{k}$ yields
$$\D^{2\lambda}_{h_{0}}\hat{\varphi}_{k,-}=-\frac{\Gamma(\mu_{k}+\frac{1}{2}+\lambda)}{\Gamma(\mu_{k}+\frac{1}{2}-\lambda)}\hat{\varphi}_{k,-}.$$
Therefore, we conclude that
$$\D_{h_{0}}^{2\lambda}=\frac{\Gamma(|\nu\cdot D_{h_{0}}|+\frac{1}{2}+\lambda)}{\Gamma(|\nu\cdot D_{h_{0}}|+\frac{1}{2}-\lambda)}\frac{\nu\cdot D_{h_{0}}}{|\nu \cdot D_{h_{0}}|}.$$
Using the fact that $\Gamma(z+1)=z\Gamma(z)$, we see that
$$\D_{h_{0}}^{2(\lambda+1)}=\Big(D_{h_{0}}^{2}-(\lambda+\frac{1}{2})^{2}\Big)\D_{h_{0}}^{2\lambda}.$$

\end{document}